\documentclass[a4paper,11pt]{amsart}
\usepackage{amsfonts,stmaryrd,amssymb,amsmath,a4wide,verbatim,hyperref,url,color,epsfig,mathrsfs,graphicx,pgf,tikz,amsthm,amsfonts,mathtext,cite,enumerate,float}
\usepackage[british]{babel}
\usetikzlibrary{arrows}
\usepackage[matrix,arrow,curve]{xy}

\makeatletter
\@addtoreset{equation}{section}\makeatother

\theoremstyle{plain}
\newtheorem{theorem}{Theorem}[section]
\theoremstyle{definition}
\newtheorem{definition}{Definition}[section]
\theoremstyle{definition}
\newtheorem{commentary}{Remark}[section]

\begin{document}

\title[Creative geometry]{Creative geometry}

\subjclass[2010]{51M04, 51N20}

\keywords{Plane geometry, Analytic geometry}

\author[Alexander Skutin]{Alexander Skutin}

\maketitle

\section{Introduction}

This work is a continuation of \cite{Sk}. As in the previous article, here we will describe some interesting ideas and a lot of new theorems in plane geometry related to them.

\bigskip

\section{Deformation of equilateral triangle}

\subsection{Deformation principle for equilateral triangle} If some triangle points lie on a circle (line) or are equivalent in the case of an equilateral triangle, then in the general case of an arbitrary triangle they are connected by some natural relations. Thus, triangle geometry can be seen as a deformation of the equilateral triangle geometry.

This principle partially describes why there exists so many relations between the Kimberling centers $X_i$. Further we will present several examples of the facts arising from this principle.

\subsection{Fermat points}

\begin{theorem}

Given a triangle $ABC$ with the first Fermat point $F$. Let $A'B'C'$ be the cevian triangle of $F$ wrt $ABC$. Let $F_A$ be the second Fermat point of $B'C'F$. Similarly define $F_B$, $F_C$. Let $M_A$ be the midpoint of $AF$. Similarly define $M_B$ and $M_C$. Then the triangles $F_AF_BF_C$ and $M_AM_BM_C$ are perspective.

\end{theorem}

\begin{theorem}

Given a triangle $ABC$ with the first Fermat point $F$. Let $F_A$ be the second Fermat point of $BCF$. Similarly define $F_B$, $F_C$. Let $F_A^*$ be the second Fermat point of $F_BF_CF$. Similarly define $F_B^*$, $F_C^*$. 
Then

\begin{enumerate}

\item $F_AF_A^*$, $F_BF_B^*$, $F_CF_C^*$ are concurrent at the first Fermat point $F'$ of $F_AF_BF_C$.

\item The second Fermat point of $F_AF_BF_C$ lies on $(F_A^*F_B^*F_C^*)$.

\item Let $F_{AB}$, $F_{BA}$, $F_{AB}'$, $F_{BA}'$ be the first Fermat points of $FF_AF_B^*$, $FF_BF_A^*$, $F'F_AF_B^*$, $F'F_BF_A^*$, respectively. Similarly define $F_{BC}$, $F_{CB}$, $F_{BC}'$, $F_{CB}'$, $F_{CA}$, $F_{AC}$, $F_{CA}'$, $F_{AC}'$. Then $F_{AB}$, $F_{AC}'$, $F_C^*$ are collinear and same for other  similar triples.

\item Let $F_{AB}F_{AC}'$ meet $F_{AB}'F_{AC}$ at $X_A$. Similarly define $X_B$, $X_C$. Then $X_A$, $X_B$, $X_C$, $F$, $F'$ are collinear.

\item{(\textbf{Miquel point})} Circles $(F_A^*X_BX_C)$, $(F_B^*X_AX_C)$, $(F_C^*X_AX_B)$, $(F_A^*F_B^*F_C^*)$ are concurrent.

\end{enumerate}

\end{theorem}

\begin{theorem}

Given a triangle $ABC$ with the first Fermat point $F_1$ and the second Fermat point $F_2$, let $A'B'C'$ be the pedal triangle of $F_1$ wrt $ABC$. Consider the isogonal conjugations $A_B$, $A_C$ of $A'$ wrt $AF_1B$, $AF_1C$, respectively. Let $F_A$ be the second Fermat point of $F_1A_BA_C$. Similarly define $F_B$, $F_C$. Then $F_2F_AF_BF_C$ is cyclic.

\end{theorem}

\subsection{Nine-point circles}

\begin{theorem}

Consider any triangle $ABC$ with the nine~- point center $N$ and centroid $G$. Let $N_A$ be the nine~- point center of $BCG$. Similarly define $N_B$, $N_C$. Then $N$ coincides with the centroid of $N_AN_BN_C$. Also $G$ lies on the nine~- point circle of $N_AN_BN_C$.

\end{theorem}

\subsection{Tangent circles}

\begin{theorem}

Given a triangle $ABC$ with the orthocenter $H$. Let the circle $\omega_A$ is tangent to $AB$, $AC$ and is externally tangent to $(BHC)$ at $A'$. Similarly define $B'$, $C'$, $\omega_B$, $\omega_C$.

\begin{enumerate}

\item Points $A'$, $B'$, $C'$, $H$ lie on the same circle.

\item Let the circle $\Omega$ is externally tangent to $\omega_A$, $\omega_B$, $\omega_C$. Then $H$ lies on $\Omega$.

\end{enumerate}

\end{theorem}

\begin{theorem}

Given a triangle $ABC$ with the incircle $\omega$. Let the circle through $B$, $C$ is tangent to $\omega$ at $T_A$. Let the circle $\pi_A$ is tangent to $\omega$ at $T_A$ and passes through $A$. Let the circle $\psi_A$ is tangent to $AB$, $AC$ and passes through $T_A$. Similarly define $\pi_B$, $\pi_C$, $\psi_B$, $\psi_C$.

\begin{enumerate}
    \item There exists a circle which is tangent to $\pi_A$, $\pi_B$, $\pi_C$, $(ABC)$.
    \item There exists a circle which is tangent to $\psi_A$, $\psi_B$, $\psi_C$, $\omega$.
\end{enumerate}

\end{theorem}

\subsection{Isogonal conjugations}

\begin{theorem}

Consider any triangle $ABC$ with the first Fermat point $F$ and let $A'B'C'$ be the cevian triangle of $F$ wrt $ABC$. Let $F_A$ be isogonal to $F$ wrt $AB'C'$ and let $F_A'$ be isogonal to $F_A$ wrt $ABC$. Similarly define $F_B'$, $F_C'$. Then the points $F_A'$, $F_B'$, $F_C'$ form an equilateral triangle with the center at $F$.

\end{theorem}

\begin{theorem}

Consider any triangle $ABC$ with the first Fermat point $F$ and let $A'B'C'$ be the cevian triangle of $F$ wrt $ABC$. Let $F_A$ be isogonal to $A$ wrt $FB'C'$ and let $F_A'$ be isogonal to $F_A$ wrt $ABC$. Let $F_A^*$ be the reflection of $F_A'$ wrt $BC$. Similarly define $F_B^*$, $F_C^*$. Then $F$ coincides with the first Fermat point of $F_A^*F_B^*F_C^*$.

\end{theorem}

\begin{theorem}

Given a triangle $ABC$ with the incenter $I$ and the incircle $\omega$ which is tangent to $BC$, $CA$, $AB$ at $A'$, $B'$, $C'$, respectively. Let $A_B$, $A_C$ be the isogonal conjugations of $A'$ wrt $AIB$, $AIC$, respectively. Similarly define $B_A$, $B_C$, $C_A$, $C_B$. Then the midpoints of $A_BA_C$, $B_AB_C$, $C_AC_B$ lie on $\omega$.

\end{theorem}

\subsection{Three equal circles}

\begin{theorem}

Consider any triangle $ABC$ with the incircle $\omega$. Let $A_1B_1C_1$ be the midpoint triangle of $ABC$ and $A'$, $B'$, $C'$ be the midpoints of smaller arcs $BC$, $CA$, $AB$ of $(ABC)$, respectively. Let $r_a$ be the radical line of the circle with diameter $A_1A'$ and $\omega$. Similarly define $r_b$, $r_c$. Then the nine~- point circle of a triangle formed by $r_a$, $r_b$, $r_c$ goes through the Feuerbach point of $ABC$.

\end{theorem}

\begin{theorem}

Consider any triangle $ABC$ with the circumcenter $O$. Consider three points $A_a$, $B_a$, $C_a$ outside wrt $ABC$ such that $BA_aC\sim AC_aB\sim CB_aA\sim BOC$. Let $A_A$, $B_A$, $C_A$ be the reflections of $A_a$, $B_a$, $C_a$ wrt $AO$, $BO$, $CO$, respectively. Let $\Omega_A$ be the circumcircle of $A_AB_AC_A$. Similarly define $\Omega_B$, $\Omega_C$. Then in the case when $\Omega_A$, $\Omega_B$, $\Omega_C$ passes through the same point we will get that they are coaxial.

\end{theorem}

\begin{commentary}
In all theorems from this section consider the case of an equilateral triangle.

\end{commentary}
\medskip

\section{Construction of midpoint analog}

\begin{definition}

For any pairs of points $A$, $B$ and $C$, $D$ denote $\mathcal{M}(AB, CD)$ as the Miquel point of the complete quadrilateral formed by the four lines $AC$, $AD$, $BC$, $BD$.

\end{definition}

Consider any two segments $AB$ and $CD$, then the point $\mathcal{M}(AB, CD)$ can be seen as midpoint between the segments $AB$, $CD$.

\begin{commentary}

In the case when $A=B$ and $C=D$ we will get that the point $\mathcal{M}(AB, CD)$ is midpoint of $AC$.

\end{commentary}

\begin{theorem}

Consider any eight lines $l_i^j$, where $1\leq i\leq 4$, $j = 1, 2$. Let given that $l_i^1\parallel l_i^2$ for any $1\leq i\leq 4$. Let $M^j$ be the Miquel point for $l_1^j$, $l_2^j$, $l_3^j$, $l_4^j$, $j = 1, 2$. By definition let $P_{pq}^j = l_p^j\cap l_q^j$ and let $Q_{pq} = \mathcal{M}(M^1M^2, P_{pq}^1P_{pq}^2)$. Then the circles $(Q_{12}Q_{23}Q_{31})$, $(Q_{12}Q_{24}Q_{41})$, $(Q_{13}Q_{34}Q_{41})$, $(Q_{23}Q_{34}Q_{42})$ passes through the same point.

\end{theorem}

\section{Combination of different facts}

Here we will build combinations of some well-known constructions from geometry.

\begin{theorem}[\textbf{Pappus and Mixtilinear circles}]

Given a triangle $ABC$, let its $A$~-- mixtilinear incircle $\omega_A$ is tangent to $(ABC)$ at $A_1$ and $A$~-- mixtilinear excircle $\Omega_A$ is tangent to $(ABC)$ at $A_2$. Similarly define $B_1$, $B_2$, $C_1$, $C_2$ and $\omega_B$, $\omega_C$, $\Omega_B$, $\Omega_C$. Let the radical line of $\omega_A$, $\omega_B$ meet $A_1B_1$ at $X$ and the radical line of $\Omega_A$, $\Omega_B$ meet $A_2B_2$ at $Y$. Let $l_C$ be the line through the points $A_1Y\cap A_2X$, $A_1B_2\cap A_2B_1$ and $XB_2\cap YB_1$. Similarly define $l_B$, $l_C$. Then the triangle formed by the lines $l_A$, $l_B$, $l_C$ is perspective to $ABC$.

\end{theorem}

\begin{theorem}[\textbf{Conics and Japanese theorem}]

Given a cyclic quadrilateral $ABCD$. Let $I_A$, $I_B$, $I_C$, $I_D$ be the incenters of $DAC$, $ABC$, $BCD$, $CDA$, respectively. Consider any point $P$ on $(ABCD)$ and two conics $\mathcal{C}_{AC}$, $\mathcal{C}_{BD}$ through $ACI_BI_DP$, $BDI_AI_CP$, respectively. Then the conics $\mathcal{C}_{AC}$, $\mathcal{C}_{BD}$ meet on $(ABCD)$ at two different points.

\end{theorem}

\begin{theorem}[\textbf{Conics and Japanese theorem 2}]

Given a cyclic quadrilateral $ABCD$. Consider the point $E = AC\cap BD$. Then the incenters of the triangles $ABC$, $BCD$, $CDA$, $DAB$, $AEB$, $BEC$, $CED$, $DEA$ lie on the same conic.

\end{theorem}

\subsection{Combinations with radical lines}

\begin{theorem}[\textbf{Japanese theorem and radical lines}]

Given a cyclic quadrilateral $ABCD$. Let $I_A$, $I_B$, $I_C$, $I_D$ be the incenters of $DAC$, $ABC$, $BCD$, $CDA$, respectively. Then the circles with diameters $AI_A$, $BI_B$, $CI_C$, $DI_D$ have the same radical center.

\end{theorem}

\begin{theorem}[\textbf{Japanese theorem and radical lines 2}]

Given a cyclic quadrilateral $ABCD$. Let $I_A$, $I_B$, $I_C$, $I_D$ be the incenters of $DAC$, $ABC$, $BCD$, $CDA$, respectively. Define the points $A_B = AB\cap I_AI_D$, $A_D = AD\cap I_AI_B$. Similarly define $B_A$, $B_C$, $C_B$, $C_D$, $D_C$, $D_A$. Then the circles $(AA_BA_D)$, $(BB_AB_C)$, $(CC_BC_D)$, $(DD_CD_A)$ have the same radical center.

\end{theorem}

\begin{theorem}[\textbf{Japanese theorem and radical lines 3}]

Given a cyclic quadrilateral $ABCD$. Let $I_A$, $I_B$, $I_C$, $I_D$ be the incenters of $DAC$, $ABC$, $BCD$, $CDA$, respectively. Define the points $A_B = AB\cap I_AI_D$, $A_D = AD\cap I_AI_B$. Similarly define $B_A$, $B_C$, $C_B$, $C_D$, $D_C$, $D_A$. Then the circles $(ABCD)$, $(I_AI_BI_CI_D)$, $(A_BB_AC_DD_C)$, $(A_DD_AB_CC_B)$ have the same radical center.

\end{theorem}

\begin{theorem}[\textbf{Gauss line theorem and radical lines}]

Consider any cyclic quadrilateral $ABCD$. Let the circle $\omega_{AB}$ is tanget to $AC$, $BD$ and is internally tangent to the smaller arc $AB$ of $(ABCD)$. Similarly define $\omega_{BC}$, $\omega_{CD}$, $\omega_{DA}$. Let $r_a$ be the radical line of $\omega_{AB}$, $\omega_{DA}$. Similarly define $r_b$, $r_c$, $r_d$. Then the circumcenter of $ABCD$ lies on the Gauss line of $r_a$, $r_b$, $r_c$, $r_d$.

\end{theorem}

\begin{theorem}[\textbf{Gauss line theorem and radical lines 2}] Consider any six points $P_1$, $P_2$, $\ldots$, $P_6$ such that $P_4P_5P_6$ lie on the same line and similarly $P_1P_4P_2$, $P_1P_5P_3$, $P_2P_3P_6$ lie on the same lines. Let $M_{25}$ be the midpoint of $P_2P_5$, similarly define $M_{61}$, $M_{43}$. Let $R_{153}$ be the radical center of $(P_1M_{25}M_{43})$, $(P_5M_{61}M_{43})$ and $(P_3M_{61}M_{25})$. Similarly define $R_{142}$, $R_{236}$, $R_{456}$. Then in the case when $M_{61}$ lies on $R_{456}R_{153}$ we get that $M_{61}=R_{456}R_{153}\cap R_{142}R_{236}$ and $M_{25} = R_{153}R_{142}\cap R_{236}R_{456}$.
\end{theorem}

\begin{theorem}[\textbf{Gauss line theorem and radical lines 3}] Consider any four points $A$, $B$, $C$, $D$. Let $AB\cap CD = E$, $BC\cap AD = F$ and let the perpendicular bisectors of the segments $AC$, $BD$, $EF$ form a triangle $PQR$. Let $L_{ABF}$ be the radical line of $(ABF)$ and $(PQR)$. Similarly define $L_{CDF}$, $L_{BCE}$, $L_{DAE}$. Consider the Gauss line $\mathcal{L}$ of $L_{ABF}$, $L_{CDF}$, $L_{BCE}$, $L_{DAE}$. Then in the limiting case when the lines $L_{ABF}$, $L_{CDF}$, $L_{BCE}$, $L_{DAE}$ are concurrent we get that $\mathcal{L}$ is parallel to the Gauss line of $AB$, $BC$, $CD$, $DA$.
\end{theorem}

\begin{theorem}[\textbf{Morley's theorem and radical lines}]

Consider any triangle $ABC$ and its three external Morley's triangles $A_AB_AC_A$, $A_BB_BC_B$, $A_CB_CC_C$ (see picture below for more details). Let $l_A$ be the radical line of $(A_AA_BA_C)$ and $(C_AB_AB_BC_C)$. Similarly define $l_B$, $l_C$. Then the lines $l_A$, $l_B$, $l_C$ form a triangle which is perspective to $ABC$.

\end{theorem}

\bigskip

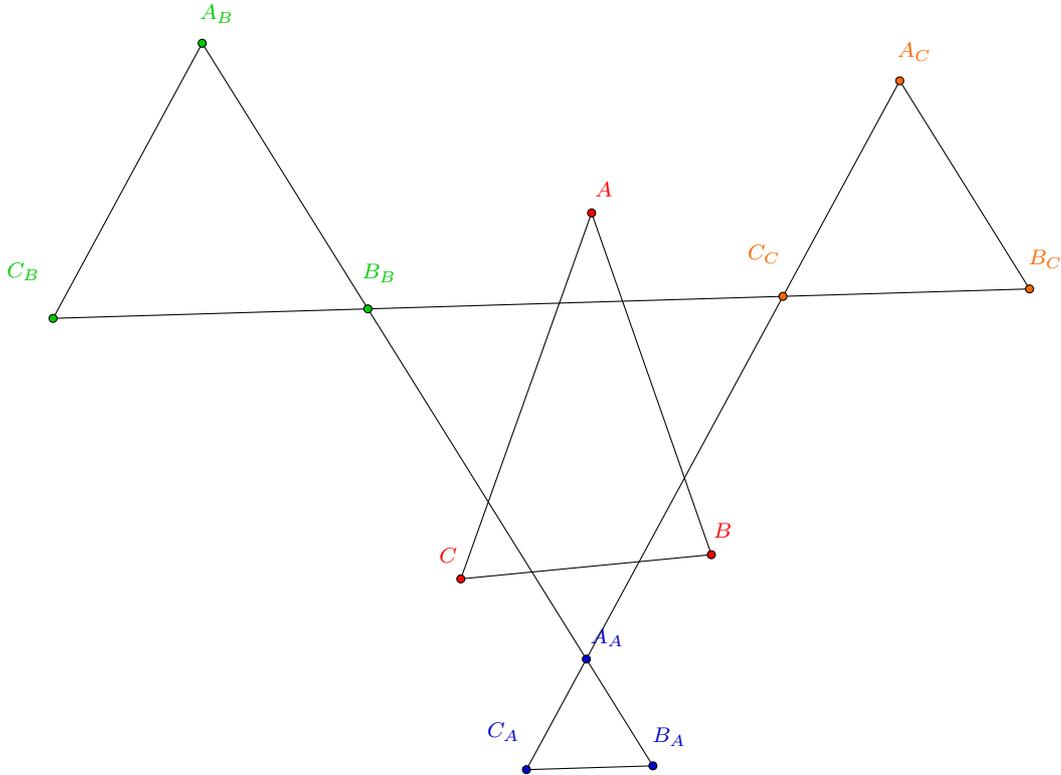
\definecolor{ffwwqq}{rgb}{1.0,0.4,0.0}
\definecolor{qqqqcc}{rgb}{0.0,0.0,0.8}
\definecolor{qqccqq}{rgb}{0.0,0.8,0.0}
\definecolor{ffqqqq}{rgb}{1.0,0.0,0.0}
\begin{tikzpicture}[line cap=round,line join=round,>=triangle 45,x=1.0cm,y=1.0cm]
\draw (6.8292322111854995,-8.603114064156882)-- (11.738886345340562,0.5309141843419108);
\draw (8.492253578166556,-8.552553554961394)-- (2.564837481332247,1.0292520311585056);
\draw (0.603813098171444,-2.619080747959214)-- (13.445980360109019,-2.228642880578118);
\draw (2.564837481332247,1.0292520311585056)-- (0.603813098171444,-2.619080747959214);
\draw (6.8292322111854995,-8.603114064156882)-- (8.492253578166556,-8.552553554961394);
\draw (11.738886345340562,0.5309141843419108)-- (13.445980360109019,-2.228642880578118);
\draw (5.965994301540906,-6.07390825740116)-- (7.685940188115136,-1.2212797216206246);
\draw (7.685940188115136,-1.2212797216206246)-- (9.26149662673033,-5.752648562408445);
\draw (5.965994301540906,-6.07390825740116)-- (9.26149662673033,-5.752648562408445);
\begin{scriptsize}
\draw [fill=ffqqqq] (7.685940188115136,-1.2212797216206246) circle (1.5pt);
\draw[color=ffqqqq] (7.844907507561563,-0.9022982946641378) node {$A$};
\draw [fill=ffqqqq] (9.26149662673033,-5.752648562408445) circle (1.5pt);
\draw[color=ffqqqq] (9.411871084962048,-5.436391434974198) node {$B$};
\draw [fill=ffqqqq] (5.965994301540906,-6.07390825740116) circle (1.5pt);
\draw[color=ffqqqq] (5.801041971821799,-5.7553728619306845) node {$C$};
\draw [fill=qqccqq] (0.603813098171444,-2.619080747959214) circle (1.5pt);
\draw[color=qqccqq] (0.21447617413311218,-1.9959489013720924) node {$C_B$};
\draw [fill=qqccqq] (2.564837481332247,1.0292520311585056) circle (1.5pt);
\draw[color=qqccqq] (2.7579532852759288,1.4217092445902657) node {$A_B$};
\draw [fill=qqccqq] (4.743874157927763,-2.4932112916572646) circle (1.5pt);
\draw[color=qqccqq] (4.892657289270794,-2.0187332890118412) node {$B_B$};
\draw [fill=qqqqcc] (8.492253578166556,-8.552553554961394) circle (1.5pt);
\draw[color=qqqqcc] (8.707872955985017,-8.170517951744085) node {$B_A$};
\draw [fill=qqqqcc] (6.8292322111854995,-8.603114064156882) circle (1.5pt);
\draw[color=qqqqcc] (6.527749717862604,-8.102164788824838) node {$C_A$};
\draw [fill=qqqqcc] (7.616956209284517,-7.137615058717248) circle (1.5pt);
\draw[color=qqqqcc] (7.890326741689114,-6.84902346863864) node {$A_A$};
\draw [fill=ffwwqq] (11.738886345340562,0.5309141843419108) circle (1.5pt);
\draw[color=ffwwqq] (11.932638579041088,0.9204527165157864) node {$A_C$};
\draw [fill=ffwwqq] (13.445980360109019,-2.228642880578118) circle (1.5pt);
\draw[color=ffwwqq] (13.658569475888,-1.8364581878938488) node {$B_C$};
\draw [fill=ffwwqq] (10.202586831312958,-2.3272511315532167) circle (1.5pt);
\draw[color=ffwwqq] (9.95690189449265,-1.7681050249746018) node {$C_C$};
\end{scriptsize}
\end{tikzpicture}

\section{Facts related to the set of confocal conics}

For any two conics $\mathcal{C}_1$, $\mathcal{C}_2$ with same focus $F$, define $\mathcal{L}_F(\mathcal{C}_1, \mathcal{C}_2)$ as in \cite[Definition 13.1]{Sk}.

\begin{theorem}

In a tangential quadrilateral $ABCD$, let $E$ be the intersection of the diagonals $AC$, $BD$. Let $I_{AB}$, $I_{BC}$, $I_{CD}$, $I_{DA}$ be the incenters of $AEB$, $BEC$, $CED$, $DEA$, respectively. Let $\mathcal{C}_{AB}$ be a conic with focuses $A$ and $B$, passing through $I_{AB}$. Similarly define the conics $\mathcal{C}_{BC}$, $\mathcal{C}_{CD}$, $\mathcal{C}_{DA}$. Let $X_{AB} = \mathcal{L}_A(\mathcal{C}_{DA}, \mathcal{C}_{AB})\cap\mathcal{L}_B(\mathcal{C}_{AB}, \mathcal{C}_{BC})$, similarly define $X_{BC}$, $X_{CD}$, $X_{DA}$. Then the lines $AB$, $CD$, $X_{BC}X_{DA}$ are concurrent.

\end{theorem}

Next fact can be seen also as a fact which is related to \cite[Section 6]{Sk}.

\begin{theorem}

Consider any triangle $ABC$ with the orthocenter $H$. Let $H_A$ be the reflection of $H$ wrt $BC$. Similarly define $H_B$ and $H_C$. Let $I_{AH_B}$ be the incenter of a triangle formed by $AC$, $H_BH_C$ and $AH_B$. And let $\mathcal{C}_{AH_B}$ be a conic with focuses $A$, $H_B$ and which goes through $I_{AH_B}$. Similarly define the conics $\mathcal{C}_{CH_B}$, $\mathcal{C}_{AH_C}$, $\mathcal{C}_{BH_C}$, $\mathcal{C}_{CH_A}$, $\mathcal{C}_{BH_A}$. Let $A'$ be the second intersection point of $\mathcal{L}_A(\mathcal{C}_{AH_B}, \mathcal{C}_{AH_C})$ with $(ABC)$. Similarly define $B'$, $C'$.

\begin{enumerate}

\item Line $\mathcal{L}_{H_A}(\mathcal{C}_{BH_A}, \mathcal{C}_{CH_A})$ is tangent to $(ABC)$ at $H_A$.

\item Lines $H_AA'$, $H_BB'$, $H_CC'$ are concurrent.
\end{enumerate}

\end{theorem}

\bigskip

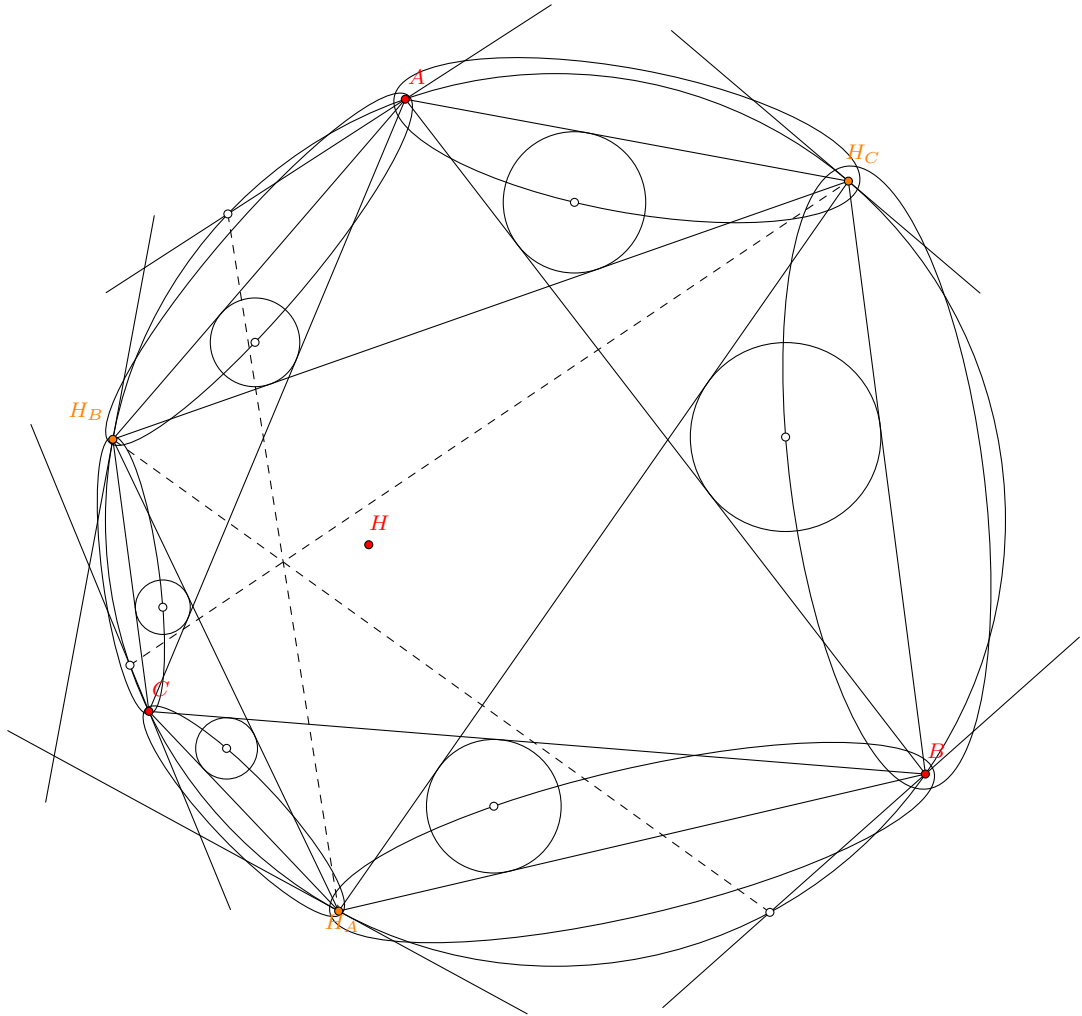
\definecolor{ffffff}{rgb}{1.0,1.0,1.0}
\definecolor{ffxfqq}{rgb}{1.0,0.4980392156862745,0.0}
\definecolor{ffqqqq}{rgb}{1.0,0.0,0.0}
\begin{tikzpicture}[line cap=round,line join=round,>=triangle 45,x=1.0cm,y=1.0cm]
\draw (3.0630285574485177,-6.705244001369529)-- (6.434903870190848,1.4125723891664603);
\draw (6.434903870190848,1.4125723891664603)-- (13.274791612467604,-7.536523702391327);
\draw (3.0630285574485177,-6.705244001369529)-- (13.274791612467604,-7.536523702391327);
\draw (2.5883360239330555,-3.0977197129704965)-- (12.262231908481686,0.325803335221328);
\draw (12.262231908481686,0.325803335221328)-- (5.558519053969256,-9.35327943816771);
\draw (5.558519053969256,-9.35327943816771)-- (2.5883360239330555,-3.0977197129704965);
\draw (6.434903870190848,1.4125723891664603)-- (12.262231908481686,0.325803335221328);
\draw (12.262231908481686,0.325803335221328)-- (13.274791612467604,-7.536523702391327);
\draw (13.274791612467604,-7.536523702391327)-- (5.558519053969256,-9.35327943816771);
\draw (5.558519053969256,-9.35327943816771)-- (3.0630285574485177,-6.705244001369529);
\draw (3.0630285574485177,-6.705244001369529)-- (2.5883360239330555,-3.0977197129704965);
\draw (2.5883360239330555,-3.0977197129704965)-- (6.434903870190848,1.4125723891664603);
\draw(8.659591817596622,0.04465836358048292) circle (0.9368671765587224cm);
\draw(4.458473613004172,-1.8116992867450217) circle (0.5884310653023586cm);
\draw(3.245598347908354,-5.324101224409755) circle (0.3611926145715426cm);
\draw(4.084330933177578,-7.196243576597974) circle (0.4065166482721754cm);
\draw(7.598720948415915,-7.963052819371638) circle (0.885655108695977cm);
\draw(11.43602797920843,-3.067554916212687) circle (1.2528745532589662cm);
\draw(8.40938834624781,-4.166755390478843) circle (5.918402438066521cm);
\draw [rotate around={-10.56403289252062:(9.34856788933628,0.8691878621938927)}] (9.34856788933628,0.8691878621938927) ellipse (3.1125939531699247cm and 0.9505452943936976cm);
\draw [rotate around={49.54105239741823:(4.511619947061953,-0.842573661902019)}] (4.511619947061953,-0.842573661902019) ellipse (3.0257324837311104cm and 0.6085657749904912cm);
\draw [rotate around={-82.50385689423618:(2.825682290690784,-4.901481857169994)}] (2.825682290690784,-4.901481857169994) ellipse (1.8572646674106268cm and 0.3735517656509554cm);
\draw [rotate around={-46.69876076748377:(4.310773805708913,-8.02926171976865)}] (4.310773805708913,-8.02926171976865) ellipse (1.872933336384172cm and 0.44495860408344623cm);
\draw [rotate around={13.248696587546782:(9.41665533321843,-8.444901570279521)}] (9.41665533321843,-8.444901570279521) ellipse (4.080455215465088cm and 0.9694064499933384cm);
\draw [rotate around={-82.66148482914447:(12.768511760474649,-3.605360183585018)}] (12.768511760474649,-3.605360183585018) ellipse (4.16247879755864cm and 1.2711663305788399cm);
\draw (9.934876600192778,2.321756301907298)-- (13.993450736659241,-1.1588996423502669);
\draw (8.353794299490254,2.663541550053717)-- (2.499849168106426,-1.1527810064983839);
\draw (3.1330028998902786,-0.13193059116345973)-- (1.704740129681054,-7.909025865268546);
\draw (1.206043055277103,-6.960860150037763)-- (8.037339410785407,-10.715809096053958);
\draw (9.820803046868967,-10.632828721828057)-- (15.299624013874432,-5.721376197145468);
\draw [dash pattern=on 3pt off 3pt] (4.1011436153702485,-0.10886002985439047)-- (5.558519053969256,-9.35327943816771);
\draw (4.1353688623889,-9.333327373691935)-- (1.5107773198892707,-2.900998606225927);
\draw [dash pattern=on 3pt off 3pt] (2.813270991202254,-6.093140078818553)-- (12.262231908481686,0.325803335221328);
\draw [dash pattern=on 3pt off 3pt] (11.229342464468905,-9.37015295265298)-- (2.5883360239330555,-3.0977197129704965);
\begin{scriptsize}
\draw [fill=ffqqqq] (3.0630285574485177,-6.705244001369529) circle (1.5pt);
\draw[color=ffqqqq] (3.220220126900938,-6.399728943478729) node {$C$};
\draw [fill=ffqqqq] (6.434903870190848,1.4125723891664603) circle (1.5pt);
\draw[color=ffqqqq] (6.583187816182196,1.7150511310985095) node {$A$};
\draw [fill=ffqqqq] (13.274791612467604,-7.536523702391327) circle (1.5pt);
\draw[color=ffqqqq] (13.415546222886523,-7.2325616353432345) node {$B$};
\draw [fill=ffqqqq] (5.953947347611348,-4.495684533636707) circle (1.5pt);
\draw[color=ffqqqq] (6.093641886729862,-4.200196449580161) node {$H$};
\draw [fill=ffxfqq] (2.5883360239330555,-3.0977197129704965) circle (1.5pt);
\draw[color=ffxfqq] (2.241128267996268,-2.7267232255121887) node {$H_B$};
\draw [fill=ffxfqq] (5.558519053969256,-9.35327943816771) circle (1.5pt);
\draw[color=ffxfqq] (5.604095957277526,-9.51751286686893) node {$H_A$};
\draw [fill=ffxfqq] (12.262231908481686,0.325803335221328) circle (1.5pt);
\draw[color=ffxfqq] (12.457738969610215,0.6900262795729637) node {$H_C$};
\draw [fill=ffffff] (8.659591817596622,0.04465836358048292) circle (1.5pt);
\draw [fill=ffffff] (4.458473613004172,-1.8116992867450217) circle (1.5pt);
\draw [fill=ffffff] (3.245598347908354,-5.324101224409755) circle (1.5pt);
\draw [fill=ffffff] (4.084330933177578,-7.196243576597974) circle (1.5pt);
\draw [fill=ffffff] (7.598720948415915,-7.963052819371638) circle (1.5pt);
\draw [fill=ffffff] (11.43602797920843,-3.067554916212687) circle (1.5pt);
\draw [fill=ffffff] (4.1011436153702485,-0.10886002985439047) circle (1.5pt);
\draw [fill=ffffff] (2.813270991202254,-6.093140078818553) circle (1.5pt);
\draw [fill=ffffff] (11.229342464468905,-9.37015295265298) circle (1.5pt);
\end{scriptsize}
\end{tikzpicture}

\section{Fun with triangle centers}

In this section we will construct some nice theorems which include some famous triangle centers. Note that the isogonal conjugation can be considered as a conjugation which is related to the incenter and the isotomic conjugation can be considered as a conjugation which is related to the centroid of $ABC$. Further, it will be useful to use the notation $P_{ABC}^{G} $ for the isotomic conjugation of the point $P$ wrt $ABC$ and $P_{ABC}^{I}$ for the isogonal conjugation of $P$ wrt $ABC$.

\begin{theorem}

Given a triangle $ABC$. Let $I$ be the incenter of $ABC$ and $G$ be the centroid of $ABC$. Let $X = II_{ABC}^{G}\cap GG_{ABC}^I$, then

\begin{enumerate}
\item Points $G$, $I$, $X_{ABC}^I$, $X_{ABC}^G$ are collinear.
\item $X_{ABC}^G$ coincides with the Spieker center of $ABC$.
\item Triangles $ABC$ and $X_{AGI}^GX_{BGI}^GX_{CGI}^G$ are perspective.
\end{enumerate}

\end{theorem}

\begin{theorem}

Given a triangle $ABC$. Let $H$ be the orthocenter of $ABC$ and $G$ be the centroid of $ABC$. Let $\triangle_{G}$ be the midpoint triangle of $ABC$ and $\triangle_{H}$ be its orthic triangle. Consider the point $X = HH_{\triangle_{G}}^{G}\cap GG_{\triangle_{H}}^{I}$. Then the point $X_{\triangle_{G}}^{G}$ coincides with the circumcenter of $ABC$.

\end{theorem}

\begin{theorem}

Given a triangle $ABC$. Let $I$ be the incenter of $ABC$ and $G$ be the centroid of $ABC$. Let $\triangle_{G}$ be the midpoint triangle of $ABC$ and $\triangle_{I}$ be the cevian triangle of $I$ wrt $ABC$. Consider the point $X = II_{\triangle_{G}}^{G}\cap GG_{\triangle_{I}}^{I}$. Then $|XG| = \frac{1}{3}\lvert XG_{\triangle_I}^I\rvert$.

\end{theorem}

\bigskip

\section{Orthogonal facts from orthogonal constructions}

Consider any situation which is related to some orthogonal objects (name  such situation as "Orthogonal"). Then one can try to add some more orthogonal objects to this situation and try to find relations between them. In the next examples we will add orthogonal hyperbolas to some well-known orthogonal situations.

\begin{theorem}[\textbf{Two orthogonal circles}]

Consider any two orthogonal circles $\omega_1$ and $\omega_2$ with centers $O_1$, $O_2$, respectively. Let the external common tangent to $\omega_1$ and $\omega_2$ is tangent to $\omega_i$ at $X_i$. Let $O_1X_2$ intersects $\omega_2$ second time at $Y_1$, similarly define $Y_2$ on $\omega_1$. Consider an orthogonal hyperbola $\mathcal{H}$ through $X_1X_2Y_1Y_2$. Let the line $X_1O_1$ meet $\mathcal{H}$ second time at $Z_1$. Similarly define $Z_2$. Then $\frac{\lvert X_1O_1\lvert}{\lvert O_1Z_1\lvert} = \frac{\lvert X_2O_2\lvert}{\lvert O_2Z_2\lvert}$.

\end{theorem}

\begin{theorem}[\textbf{Three orthogonal circles}]

Consider any three pairwise orthogonal circles $\omega_1$, $\omega_2$, $\omega_3$ with centers $O_1$, $O_2$, $O_3$, respectively. Let $A_1$, $A_1'$ be two intersection points of $\omega_2$, $\omega_3$. Similarly define $A_2$, $A_2'$, $A_3$, $A_3'$. Consider an orthogonal hyperbola $\mathcal{H}_1$ through $A_2A_2'A_3A_3'$. Similarly define $\mathcal{H}_2$ and $\mathcal{H}_3$. Let $B_1$, $B_1'$ be two intersection points of $\mathcal{H}_2$ and $\mathcal{H}_3$ different from $A_1$, $A_1'$. Similarly define $B_2$, $B_2'$, $B_3$, $B_3'$. Then

\begin{enumerate}

\item Lines $A_1A_1'$, $B_2B_2'$, $B_3B_3'$ are concurrent.

\item Points $B_1$, $B_1'$, $B_2$, $B_2'$, $B_3$, $B_3'$ lie on the same circle.

\item Lines $B_1B_1'$ and $A_1A_1'$ are orthogonal.

\end{enumerate}

\end{theorem}

\begin{theorem}[\textbf{Orthocenter construction}]

Consider any triangle with the orthocenter $H$ and altitudes $AA'$, $BB'$, $CC'$. Let $\mathcal{H}_A$ be an orthogonal hyperbola through $AHB'C'$. Similarly define $\mathcal{H}_B$, $\mathcal{H}_C$. Let $A_1B_1C_1$ be the midpoint triangle of $ABC$. Then the lines $A_1B_1$, $B_1C_1$, $C_1A_1$ pass through some intersection points of $\mathcal{H}_A$, $\mathcal{H}_B$, $\mathcal{H}_C$.

\end{theorem}

\begin{theorem}[\textbf{Perpendiculars through the orthocenter}]

Consider any triangle with the orthocenter $H$ and let $P$ be any point on the plane. Let $\ell_1$, $\ell_2$ be two perpendicular lines through $H$. Consider the intersection points $A_1$, $B_1$, $C_1$ of $\ell_1$ with $BC$, $CA$, $AB$, respectively. Similarly define $A_2$, $B_2$, $C_2$. Let $\mathcal{H}_A$ be an orthogonal hyperbola through $A_1A_2HP$. Similarly define $\mathcal{H}_B$, $\mathcal{H}_C$. Then $\mathcal{H}_A$, $\mathcal{H}_B$, $\mathcal{H}_C$ intersect at four different points.

\end{theorem}

\begin{theorem}[\textbf{Two parabolas with orthogonal directrices}]

Consider any two parabolas $\mathcal{P}_1$, $\mathcal{P}_2$ with orthogonal directrices $d_1$, $d_2$, respectively. Let $F_1$, $F_2$ be the focuses of $\mathcal{P}_1$, $\mathcal{P}_2$, respectively. Consider the four intersection points $A$, $B$, $C$, $D$ of $\mathcal{P}_1$, $\mathcal{P}_2$. Let $\mathcal{H}$ be an orthogonal hyperbola through $ABCD$. Let $P_i$ be a pole of $d_i$ wrt $\mathcal{H}$. Then $P_1P_2\perp F_1F_2$.

\end{theorem}

\begin{theorem}[\textbf{Square}]

Consider a square $ABCD$ and any point $P$ inside it. Let $l_1$, $l_2$ be two orthogonal lines through $P$. Let $X = AB\cap l_1$, $Y = CD\cap l_1$, $Z = AD\cap l_2$, $T = BC\cap l_2$. Consider an orthogonal hyperbola $\mathcal{H}_{PAXZ}$ through $PAXZ$, similarly define $\mathcal{H}_{PCYT}$. Then the internal angle bisector of $\angle XPZ$ passes through the intersection point of $\mathcal{H}_{PAXZ}$ and $\mathcal{H}_{PCYT}$ different from $P$.

\end{theorem}

\bigskip

\section{Combination of several ideas}

The next fact shows how it is possible to combine the following three ideas :

\begin{enumerate}

\item Combination of Gauss line theorem with radical lines.

\item Midpoint analog.

\item Blow-ups with circles.

\end{enumerate}

\begin{theorem}

Consider any six circles $\omega_1$, $\omega_2$, $\ldots$, $\omega_6$. And let $O_i$ be the center of $\omega_i$. Let given that $\omega_4$, $\omega_5$, $\omega_6$ share the same two external tangents. Also let given a similar property for triples of circles $\omega_1$, $\omega_4$, $\omega_2$; $\omega_1$, $\omega_5$, $\omega_3$ and $\omega_2$, $\omega_3$, $\omega_6$. Let $T_1^1$, $T_1^2$ be the common external tangents to $\omega_1$, $\omega_2$ and $L_1^1$, $L_1^2$ be the common external tangents to $\omega_1$, $\omega_3$. Consider $X_1 = T_1^1\cap L_1^1$ and $Y_1 = T_1^2\cap L_1^2$. Similarly define $X_i$, $Y_i$ (see picture below for more details). Let $M_{25} = \mathcal{M}(X_2Y_2, X_5Y_5)$, similarly define $M_{61}$, $M_{43}$ (see definition of $\mathcal{M}$ in Section 3). Let $R_{153}$ be the radical center of $(O_1M_{25}M_{43})$, $(O_5M_{61}M_{43})$ and $(O_3M_{61}M_{25})$. Similarly define $R_{142}$, $R_{236}$, $R_{456}$. Then in the case when $M_{61}$ lies on $R_{456}R_{153}$, we get that $M_{61} = R_{456}R_{153}\cap R_{142}R_{236}$ and $M_{25} = R_{153}R_{142}\cap R_{236}R_{456}$.

\end{theorem}
\bigskip
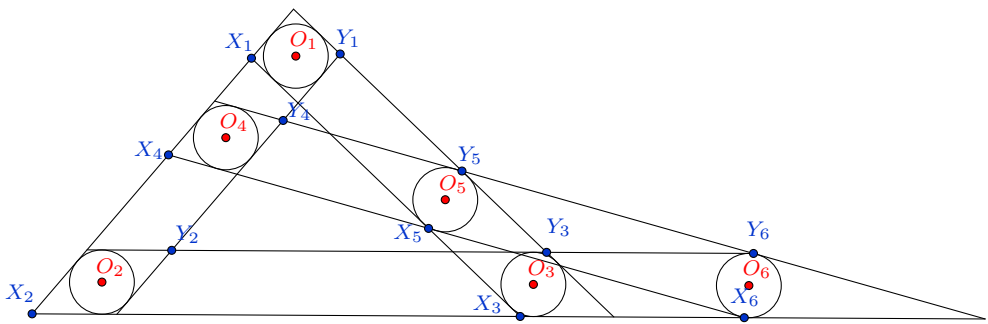
\definecolor{qqttcc}{rgb}{0.0,0.2,0.8}
\definecolor{ffqqqq}{rgb}{1.0,0.0,0.0}
\begin{tikzpicture}[line cap=round,line join=round,>=triangle 45,x=1.0cm,y=1.0cm]
\draw(9.897059284384463,0.5938594366963297) circle (0.4261613726502333cm);
\draw(7.347003120571587,-2.4028673684625046) circle (0.4261613726502333cm);
\draw(13.021534556236974,-2.4336339867140584) circle (0.4261613726502333cm);
\draw(15.855357682306606,-2.4489986305266247) circle (0.4261613726502186cm);
\draw(11.862555178477496,-1.3106286374175333) circle (0.42616137265022325cm);
\draw(8.976601090728739,-0.4878272261698926) circle (0.42616137265026854cm);
\draw (6.4290193739943104,-2.8240578088912494)-- (18.96355156117395,-2.892018515627479);
\draw (18.96355156117395,-2.892018515627479)-- (8.829747747604987,-0.0028152138468171242);
\draw (6.4290193739943104,-2.8240578088912494)-- (9.867123997630712,1.2162690006902583);
\draw (9.867123997630713,1.2162690006902577)-- (14.07968724664235,-2.865538798198623);
\draw (7.1509807177576405,-1.9756369223845216)-- (15.916053685487714,-2.0231600803760763);
\draw (10.480511147904977,0.6219210216266451)-- (7.543025523385536,-2.8300978145404887);
\draw (8.222812942294691,-0.7160613869576476)-- (15.794661679125504,-2.8748371806771753);
\draw (9.313607420863953,0.5657978517660147)-- (12.847970369344837,-2.8588605793205626);
\begin{scriptsize}
\draw [fill=ffqqqq] (9.897059284384463,0.5938594366963297) circle (1.5pt);
\draw[color=ffqqqq] (10.005033553457302,0.8047526036475704) node {$O_1$};
\draw [fill=ffqqqq] (8.976601090728739,-0.4878272261698926) circle (1.5pt);
\draw[color=ffqqqq] (9.083031986137785,-0.28626222796774314) node {$O_4$};
\draw [fill=ffqqqq] (11.862555178477496,-1.3106286374175333) circle (1.5pt);
\draw[color=ffqqqq] (11.970352683796273,-1.1105845451882024) node {$O_5$};
\draw [fill=ffqqqq] (7.347003120571587,-2.4028673684625046) circle (1.5pt);
\draw[color=ffqqqq] (7.457397643758635,-2.201599376803516) node {$O_2$};
\draw [fill=ffqqqq] (13.021534556236974,-2.4336339867140584) circle (1.5pt);
\draw[color=ffqqqq] (13.134986242515662,-2.225844150839412) node {$O_3$};
\draw [fill=ffqqqq] (15.855357682306606,-2.4489986305266247) circle (1.5pt);
\draw[color=ffqqqq] (15.961648942324183,-2.2379665378573597) node {$O_6$};
\draw [fill=qqttcc] (9.313607420863953,0.5657978517660147) circle (1.5pt);
\draw[color=qqttcc] (9.155821583557746,0.7805078296116746) node {$X_1$};
\draw [fill=qqttcc] (10.480511147904977,0.6219210216266451) circle (1.5pt);
\draw[color=qqttcc] (10.587350332816996,0.8289973776834663) node {$Y_1$};
\draw [fill=qqttcc] (8.222812942294691,-0.7160613869576476) circle (1.5pt);
\draw[color=qqttcc] (7.966924825698368,-0.6499338385061811) node {$X_4$};
\draw [fill=qqttcc] (9.730389239162736,-0.2595930653821998) circle (1.5pt);
\draw[color=qqttcc] (9.944375555607333,-0.14079358375236803) node {$Y_4$};
\draw [fill=qqttcc] (11.642830066394666,-1.691127112412626) circle (1.5pt);
\draw[color=qqttcc] (11.375904304866584,-1.7530710571394426) node {$X_5$};
\draw [fill=qqttcc] (12.082280290560396,-0.9301301624225062) circle (1.5pt);
\draw[color=qqttcc] (12.188721476056157,-0.7226681606138687) node {$Y_5$};
\draw [fill=qqttcc] (6.4290193739943104,-2.8240578088912494) circle (1.5pt);
\draw[color=qqttcc] (6.2685008858992575,-2.553148600324006) node {$X_2$};
\draw [fill=qqttcc] (8.264986867148867,-1.9816769280337607) circle (1.5pt);
\draw[color=qqttcc] (8.476452007638102,-1.7651934441573904) node {$Y_2$};
\draw [fill=qqttcc] (12.847970369344837,-2.8588605793205626) circle (1.5pt);
\draw[color=qqttcc] (12.419221867886037,-2.6986172445393812) node {$X_3$};
\draw [fill=qqttcc] (13.195098743129117,-2.008407394107556) circle (1.5pt);
\draw[color=qqttcc] (13.353355034775548,-1.680336735031755) node {$Y_3$};
\draw [fill=qqttcc] (15.794661679125504,-2.8748371806771753) circle (1.5pt);
\draw[color=qqttcc] (15.803938147914264,-2.6016381483957978) node {$X_6$};
\draw [fill=qqttcc] (15.916053685487714,-2.0231600803760763) circle (1.5pt);
\draw[color=qqttcc] (15.98591214146417,-1.716703896085599) node {$Y_6$};
\end{scriptsize}
\end{tikzpicture}

\bigskip

The next fact shows how it is possible to combine the following three ideas :

\begin{enumerate}

\item Deformation of equilateral triangle.

\item Midpoint analog.

\item Combination of Morley's theorem with radical lines.

\end{enumerate}

\begin{theorem}

Consider any triangle $ABC$ and its internal and external Morley's triangles $A_1B_1C_1$, $A_2B_2C_2$, respectively (see picture below for more details). Let $M_A = \mathcal{M}(B_1B_2, C_1C_2)$ (see definition of $\mathcal{M}$ in Section 3). Similarly define $M_B$, $M_C$. Let $R_A$ be the radical line of $(B_2C_2M_A)$ and $(A_2M_BM_C)$. Similarly define $R_B$, $R_C$. Let $Q_A = B_1C_1\cap R_A$, similarly define $Q_B$, $Q_C$. Then the triangles $A_1B_1C_1$ and $Q_AQ_BQ_C$ are perspective at a point on $(A_1B_1C_1)$.

\end{theorem}

\bigskip

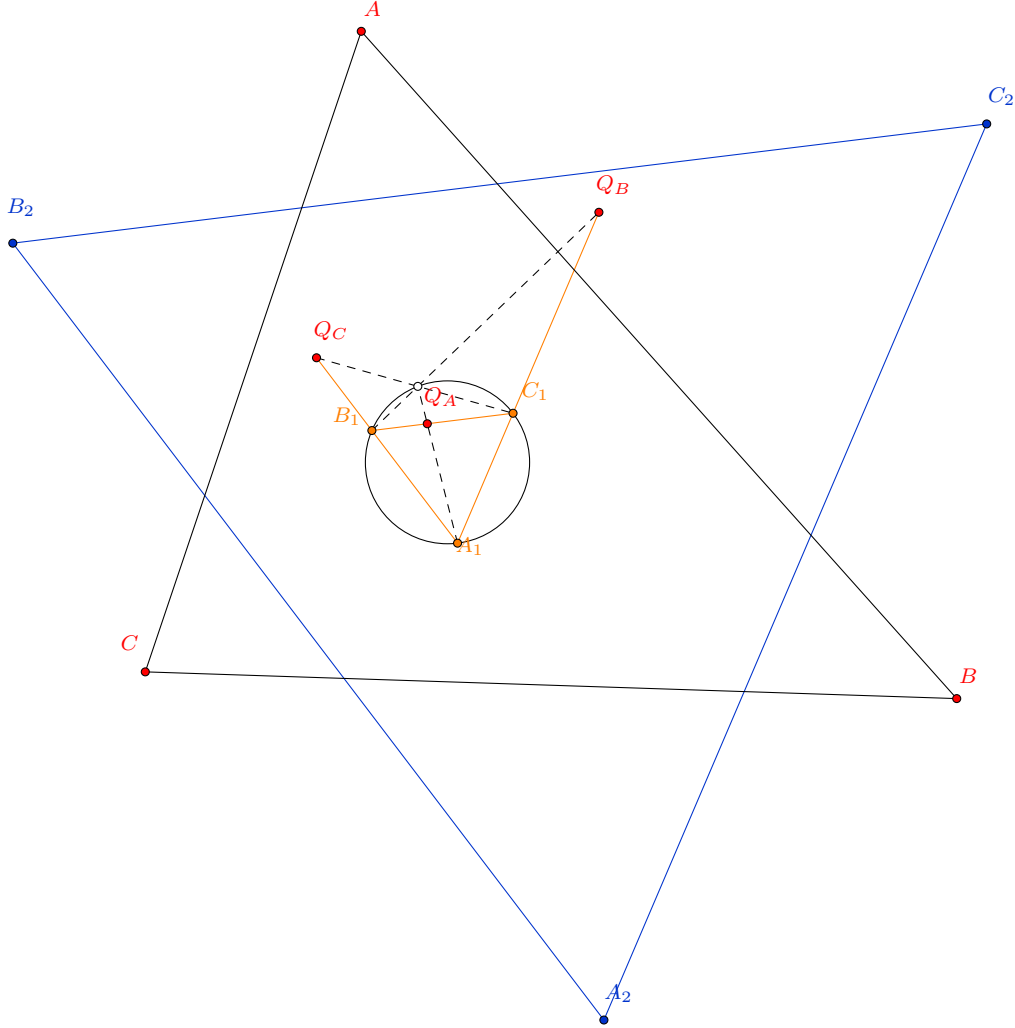
\definecolor{ffffff}{rgb}{1.0,1.0,1.0}
\definecolor{qqttcc}{rgb}{0.0,0.2,0.8}
\definecolor{ffxfqq}{rgb}{1.0,0.4980392156862745,0.0}
\definecolor{ffqqqq}{rgb}{1.0,0.0,0.0}
\begin{tikzpicture}[line cap=round,line join=round,>=triangle 45,x=1.0cm,y=1.0cm]
\draw [color=ffxfqq] (9.58,-0.48)-- (11.436524523145138,-0.2507679932037383);
\draw [color=ffxfqq] (11.436524523145138,-0.2507679932037383)-- (10.70678300281862,-1.9731813963943496);
\draw [color=qqttcc] (12.631347328423733,-8.296612815417967)-- (17.6657391560871,3.5860932104101018);
\draw [color=qqttcc] (4.857817958185923,2.0046514128574286)-- (17.6657391560871,3.5860932104101018);
\draw [color=qqttcc] (12.631347328423733,-8.296612815417967)-- (4.857817958185923,2.0046514128574286);
\draw(10.574435841987919,-0.9013164631993626) circle (1.0800047249864428cm);
\draw [color=ffxfqq] (11.436524523145138,-0.2507679932037383)-- (12.565632626368538,2.414272814359117);
\draw [dash pattern=on 3pt off 3pt] (9.58,-0.48)-- (12.565632626368538,2.414272814359117);
\draw [dash pattern=on 3pt off 3pt] (8.851761189908432,0.48504175217499684)-- (11.436524523145138,-0.2507679932037383);
\draw (6.6,-3.68)-- (17.27156480382235,-4.034443310159925);
\draw (17.27156480382235,-4.034443310159925)-- (9.440498984580536,4.8139617297374295);
\draw (9.440498984580536,4.8139617297374295)-- (6.6,-3.68);
\draw [color=ffxfqq] (10.70678300281862,-1.9731813963943496)-- (8.851761189908432,0.48504175217499684);
\draw [dash pattern=on 3pt off 3pt] (10.70678300281862,-1.9731813963943496)-- (10.18421765592714,0.10572870622476381);
\begin{scriptsize}
\draw [fill=ffqqqq] (6.6,-3.68) circle (1.5pt);
\draw[color=ffqqqq] (6.394224666267456,-3.2954305934960155) node {$C$};
\draw [fill=ffqqqq] (9.440498984580536,4.8139617297374295) circle (1.5pt);
\draw[color=ffqqqq] (9.58756050339565,5.108084767367658) node {$A$};
\draw [fill=ffxfqq] (9.58,-0.48) circle (1.5pt);
\draw[color=ffxfqq] (9.251419888961104,-0.2911738519872524) node {$B_1$};
\draw [fill=ffqqqq] (17.27156480382235,-4.034443310159925) circle (1.5pt);
\draw[color=ffqqqq] (17.423838577401025,-3.7366151499413585) node {$B$};
\draw [fill=ffxfqq] (10.70678300281862,-1.9731813963943496) circle (1.5pt);
\draw[color=ffxfqq] (10.869096595927362,-2.0138945009643057) node {$A_1$};
\draw [fill=ffxfqq] (11.436524523145138,-0.2507679932037383) circle (1.5pt);
\draw[color=ffxfqq] (11.730456920415888,0.044966762447294584) node {$C_1$};
\draw [fill=qqttcc] (4.857817958185923,2.0046514128574286) circle (1.5pt);
\draw[color=qqttcc] (4.965627054920632,2.48198621709776) node {$B_2$};
\draw [fill=qqttcc] (12.631347328423733,-8.296612815417967) circle (1.5pt);
\draw[color=qqttcc] (12.822913917328165,-7.938372830373195) node {$A_2$};
\draw [fill=qqttcc] (17.6657391560871,3.5860932104101018) circle (1.5pt);
\draw[color=qqttcc] (17.865023133846364,3.952601405248903) node {$C_2$};
\draw [fill=ffqqqq] (8.851761189908432,0.48504175217499684) circle (1.5pt);
\draw[color=ffqqqq] (9.041332004939513,0.8433007217293436) node {$Q_C$};
\draw [fill=ffqqqq] (10.30882873949559,-0.3900087973616805) circle (1.5pt);
\draw[color=ffqqqq] (10.490938404688496,-0.03906839116134214) node {$Q_A$};
\draw [fill=ffqqqq] (12.565632626368538,2.414272814359117) circle (1.5pt);
\draw[color=ffqqqq] (12.759887552121686,2.7761092547279884) node {$Q_B$};
\draw [fill=ffffff] (10.18421765592714,0.10572870622476381) circle (1.5pt);
\end{scriptsize}
\end{tikzpicture}

\bigskip

\begin{commentary}

In the previous theorem, consider the case of an equilateral triangle.

\end{commentary}

\bigskip
\section{Some new line and its properties}

\subsection{Construction} Consider any two (oriented) segments $A_1A_2$ and $B_1B_2$, let the circle $\omega_A$ is tangent to $A_1A_2$, $B_1A_1$ and $A_2B_2$. Also let the circle $\omega_B$ is tangent to $A_1B_1$, $B_1B_2$, $B_2A_2$. Consider the radical line $R_{A_1A_2, B_1B_2}$ of $\omega_A$, $\omega_B$. Next we will introduce some nice-looking situations with different pairs of (oriented) segments $i_{\alpha}$, $j_{\alpha}$ and lines $R_{i_{\alpha}, j_{\alpha}}$.

\begin{theorem}

Given a triangle $ABC$ with the circumcircle $(ABC)$. Let $AA'$, $BB'$, $CC'$ be the diameters of $(ABC)$. Then the lines $R_{AC', A'B}$, $R_{CB', A'B}$, $R_{AC', B'C}$ are concurrent.

\end{theorem}

\bigskip

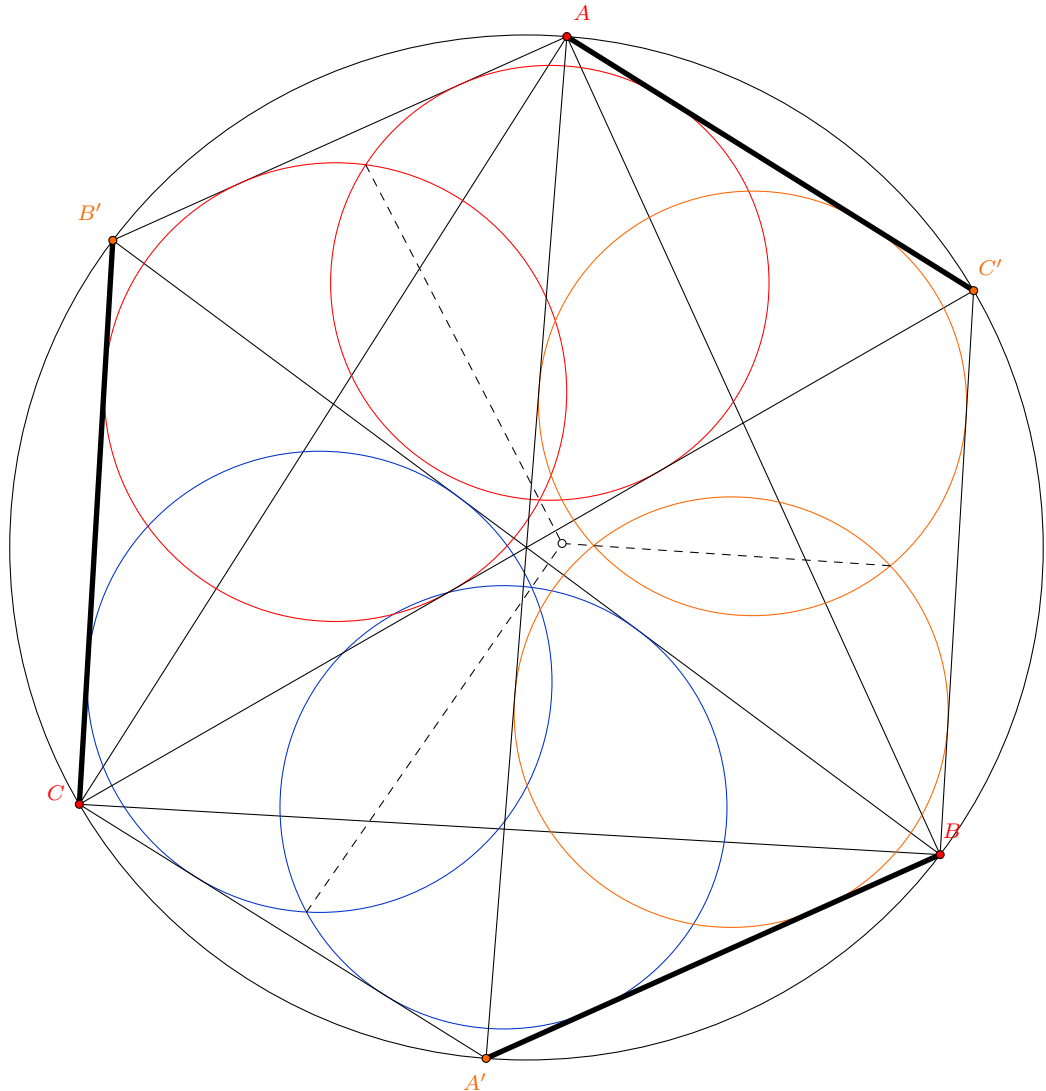
\definecolor{ffffff}{rgb}{1.0,1.0,1.0}
\definecolor{qqttcc}{rgb}{0.0,0.2,0.8}
\definecolor{ffwwqq}{rgb}{1.0,0.4,0.0}
\definecolor{ffqqqq}{rgb}{1.0,0.0,0.0}
\begin{tikzpicture}[line cap=round,line join=round,>=triangle 45,x=1.0cm,y=1.0cm]
\draw (-11.459366439909262,1.7188173940294336)-- (-5.047375122169516,11.899878409712272);
\draw (-5.047375122169516,11.899878409712272)-- (-0.13832942176329324,1.0515440454881493);
\draw (-11.459366439909262,1.7188173940294336)-- (-0.13832942176329324,1.0515440454881493);
\draw (-11.018546030768688,9.197827579502086)-- (-5.047375122169516,11.899878409712272);
\draw(-5.5784377262659905,5.124685812495118) circle (6.795973971319774cm);
\draw (-11.018546030768688,9.197827579502086)-- (-0.13832942176329324,1.0515440454881493);
\draw [color=qqttcc] (-5.880895543554122,1.6792187431791097) circle (2.939338447929837cm);
\draw [color=qqttcc] (-8.300612933592214,3.343042084349495) circle (3.0577133554244558cm);
\draw [color=ffqqqq] (-8.090521741126603,7.184487023006913) circle (3.041414193354866cm);
\draw [color=ffqqqq] (-5.27159837784003,8.633971738295285) circle (2.88300190697625cm);
\draw (-11.459366439909262,1.7188173940294336)-- (0.3024909873772823,8.530554230960803);
\draw [color=ffwwqq] (-2.6052789403691037,7.035384791189447) circle (2.814758076181243cm);
\draw [color=ffwwqq] (-2.8863731849397642,2.9411280830929445) circle (2.854463788040125cm);
\draw (-6.109500330362463,-1.650506784722037)-- (-5.047375122169516,11.899878409712272);
\draw [dash pattern=on 3pt off 3pt] (-8.471271065266354,0.2900948581376538)-- (-5.110242695457451,5.17807700960643);
\draw [dash pattern=on 3pt off 3pt] (-0.7929908867953885,4.88167287027209)-- (-5.110242695457451,5.17807700960643);
\draw [dash pattern=on 3pt off 3pt] (-7.6923649497348325,10.199726891576091)-- (-5.110242695457451,5.17807700960643);
\draw [line width=2.0pt] (-5.047375122169516,11.899878409712272)-- (0.3024909873772823,8.530554230960803);
\draw (0.3024909873772823,8.530554230960803)-- (-0.13832942176329324,1.0515440454881493);
\draw [line width=2.0pt] (-0.13832942176329324,1.0515440454881493)-- (-6.109500330362463,-1.650506784722037);
\draw (-6.109500330362463,-1.650506784722037)-- (-11.459366439909262,1.7188173940294336);
\draw [line width=2.0pt] (-11.459366439909262,1.7188173940294336)-- (-11.018546030768688,9.197827579502086);
\begin{scriptsize}
\draw [fill=ffqqqq] (-11.459366439909262,1.7188173940294336) circle (1.5pt);
\draw[color=ffqqqq] (-11.76384219408894,1.8774460444105618) node {$C$};
\draw [fill=ffqqqq] (-5.047375122169516,11.899878409712272) circle (1.5pt);
\draw[color=ffqqqq] (-4.846550525334716,12.212381868763995) node {$A$};
\draw [fill=ffqqqq] (-0.13832942176329324,1.0515440454881493) circle (1.5pt);
\draw[color=ffqqqq] (0.017867486885995675,1.3640475045398728) node {$B$};
\draw [fill=ffwwqq] (-6.109500330362463,-1.650506784722037) circle (1.5pt);
\draw[color=ffwwqq] (-6.252322703178317,-1.9618821667963293) node {$A'$};
\draw [fill=ffwwqq] (-11.018546030768688,9.197827579502086) circle (1.5pt);
\draw[color=ffwwqq] (-11.317565312233828,9.578424142470896) node {$B'$};
\draw [fill=ffwwqq] (0.3024909873772823,8.530554230960803) circle (1.5pt);
\draw[color=ffwwqq] (0.5310859010193736,8.84180884613469) node {$C'$};
\draw [fill=ffffff] (-5.110242695457451,5.17807700960643) circle (1.5pt);
\end{scriptsize}
\end{tikzpicture}

\begin{theorem}

Consider any triangle $ABC$. Let the incircle of $ABC$ is tangent to its sides at $A'$, $B'$, $C'$. Then the lines $R_{AB', BA'}$, $R_{AC', CA'}$, $R_{BC', CB'}$ are concurrent. 

\end{theorem}

\begin{theorem}

Consider any two segments $PP'$, $QQ'$ and let $X = PQ'\cap P'Q$. Consider the case when $X$ lies on $R_{PP', QQ'}$. Let $Y$ be any other point on $R_{PP', QQ'}$ and $R' = PY\cap P'Q'$, $R = P'Y\cap PQ$. Then $R_{PP', QQ'} = R_{PP', RR'}$.

\end{theorem}

\begin{theorem}

Given a cyclic quadrilateral $ABCD$. Let $W_{AB}$, $W_{BC}$, $W_{CD}$, $W_{DA}$ be the midpoints of smaller arcs $AB$, $BC$, $CD$, $DA$ of $(ABCD)$. Then the circumcenter of $ABCD$, point $R_{AB, CD}\cap R_{BC, DA}$ and the point $W_{AB}W_{CD}\cap W_{BC}W_{DA}$ are collinear.

\end{theorem}

\bigskip

\section{On the symmedian points in the incircle construction}

The main idea of this section is to construct a lot of nice facts which include the incenter and the symmedian points.

\begin{theorem}

Given a triangle $ABC$ with the incenter $I$ and let $S_A$, $S_B$, $S_C$ be the symmedian points of $BCI$, $CAI$, $ABI$, respectively. Let the lines $AI$, $BI$, $CI$ meet the sides of $S_AS_BS_C$ at $A'$, $B'$, $C'$, respectively. Then $S_AA'$, $S_BB'$, $S_CC'$ are concurrent on $SI$ where $S$ is the symmedian point of $ABC$.

\end{theorem}

\begin{theorem}

Given a triangle $ABC$ with the incenter $I$. Let $A'B'C'$ be the circumcevian triangle of $I$ wrt $ABC$. Let $S_{AB}$, $S_{AC}$ be the symmedian points of $A'IB$, $A'IC$, respectively. Let $l_A$ be the $A'$~-- median of $A'S_{AB}S_{AC}$, similarly define $l_B$, $l_C$. Then $l_A$, $l_B$, $l_C$ are concurrent.

\end{theorem}

\begin{theorem}

Consider any triangle $ABC$ with the incenter $I$ and let $A'B'C'$ be the cevian triangle of $I$ wrt $ABC$.

\begin{enumerate}
    \item Let the line $\ell_A$ goes through the symmedian points of $AA'B$ and $AA'C$. Similarly define $\ell_B$, $\ell_C$. Then $\ell_A$, $\ell_B$, $\ell_C$ form a triangle which is perspective to $ABC$.
    \item Let the line $\mu_A$ goes through the symmedian points of $IA'B$ and $IA'C$. Similarly define $\mu_B$, $\mu_C$. Then $\mu_A$, $\mu_B$, $\mu_C$ form a triangle which is perspective to $XYZ$, where $X$, $Y$, $Z$ are the midpoints of $IA$, $IB$, $IC$, respectively.
\end{enumerate}

\end{theorem}

\begin{theorem}

Consider any triangle $ABC$ with the incenter $I$ and let $A'B'C'$ be the circumcevian triangle of $I$ wrt $ABC$. Let $B'C'$ meet $AB$ and $AC$ at $A^B$, $A^C$, respectively. Similarly define $B^A$, $B^C$, $C^A$, $C^B$.

\begin{enumerate}
    \item Let $S_A$ be the symmedian point of $IB^CC^B$. Similarly define $S_B$, $S_C$. Then $S_AS_BS_C$ is perspective to both $ABC$ and $A'B'C'$.
    \item Consider the intersection point $P^A$ of $(ABC)$, $(A'B^CC^B)$ different from $A'$. Let $\ell_A$ be the line through the symmedian points of $A'B^CC^B$, $P^AB^CC^B$. Similarly define $\ell_B$, $\ell_C$. Then $\ell_A$, $\ell_B$, $\ell_C$ form a triangle which is perspective to $ABC$.
    \item Let $\eta_A$ be the line through the symmedian points of $A'P^AB^C$, $A'P^AC^B$. Similarly define $\eta_B$, $\eta_C$. Then $\eta_A$, $\eta_B$, $\eta_C$ form a triangle which is perspective to $ABC$.
\end{enumerate}

\end{theorem}

\begin{theorem}

Let $I$ be the incenter of $ABC$. Let the perpendicular line to $AI$ through $I$ intersect $AB$, $AC$ at $A_B$, $A_C$, respectively. Similarly define $B_A$, $B_C$, $C_A$, $C_B$.

\begin{enumerate}
    \item Let the line $L_A$ goes through the symmedian points of $AC_AB_A$, $AA_BA_C$. Similarly define $L_B$, $L_C$. Then $L_A$, $L_B$, $L_C$ are concurrent.
    \item Let $S_A$ and $S_A'$ be the symmedian points of $IC_BB_C$ and $AC_AB_A$, respectively. Similarly define $S_B$, $S_C$, $S_B'$, $S_C'$. Then the perpendiculars from $S_A$, $S_B$, $S_C$ to $S_B'S_C'$, $S_A'S_C'$, $S_A'S_B'$ are concurrent.
\end{enumerate}

\end{theorem}

\begin{theorem}

Given a triangle $ABC$ and let $I$, $I_A$, $I_B$, $I_C$ be the incenter and $A$, $B$, $C$~-- excenters of $ABC$, respectively. Consider any point $P$ and let $S_{PXY}$ be the symmedian point of $PXY$, for arbitrary points $X$, $Y$. Then the triangles $S_{PIA}S_{PIB}S_{PIC}$ and $S_{PII_A}S_{PII_B}S_{PII_C}$ are perspective.

\end{theorem}

\begin{theorem}[\textbf{Quadrilateral version}]

Consider any tangental quadrilateral $ABCD$ which has circumcircle $\omega$. Let $P = AC\cap BD$ and $I_{AB}$, $I_{BC}$, $I_{CD}$, $I_{DA}$ be the incenters of $PAB$, $PBC$, $PCD$, $PDA$. Let $\ell_A$ be the line through the symmedian points of $AI_{AB}I_{DA}$, $AI_{CD}I_{BC}$. Similarly define $\ell_B$, $\ell_C$, $\ell_D$. Then $P$, $\ell_A\cap\ell_C$, $\ell_B\cap\ell_D$ are collinear.

\end{theorem}

\section{Orthocenter, Circumcenter and Isogonal conjugate points}

It is well-known that the orthocenter is isogonal conjugated to the circumcenter of a triangle. Here we will introduce more sophisticated facts which are related to the orthocenter~-- circumcenter pair and to the isogonal conjugation.

\begin{theorem}

Given a triangle $ABC$ with the orthocenter $H$ and the circumcenter $O$. Let $P$, $Q$ be some isogonal conjugated points wrt $ABC$. Let $H_A$, $H_B$, $H_C$ be the isogonal conjugations of $H$ wrt $BCP$, $ACP$, $ABP$, respectively. Let $O_A$, $O_B$, $O_C$ be the isogonal conjugations of $O$ wrt $BCQ$, $ACQ$, $ABQ$, respectively. Then the circles $(H_AH_BH_C)$, $(O_AO_BO_C)$ have the same radiuses.

\end{theorem}

\begin{theorem}

Consider any triangle $ABC$ with the orthocenter $H$ and the circumcenter $O$. Let $\mathcal{C}$ be the conic with focuses $H$, $O$ and which is tangent to the sides of $ABC$. Let the perpendicular bisector of $OH$ meet $\mathcal{C}$ at $X$, $Y$. Consider the isogonal conjugations $X'$, $Y'$ of $X$, $Y$ wrt $ABC$, respectively. Then

\begin{enumerate}
    \item The lines $XY'$, $X'Y$ meet on $(ABC)$.
    \item $XY\parallel X'Y'$
\end{enumerate}

\end{theorem}

\begin{definition}

For any two segments $AB$, $CD$, let $KH(AB, CD)$ be the Kantor~- Hervey point of the complete quadrilateral formed by the lines $AC$, $AD$, $BC$, $BD$.

\end{definition}

\begin{theorem}

Consider any quadrilateral $ABCD$ and let $E = AD\cap BC$, $F = AB\cap CD$. Consider a conic $\mathcal{C}$ with focuses $P$, $Q$ which is inscribed in $ABCD$. Then the points $KH(PQ, AC)$, $KH(PQ, BD)$, $KH(PQ, EF)$ are collinear.

\end{theorem}

\begin{theorem}

Consider any triangle $ABC$ with the orthocenter $H$ and the circumcenter $O$. Let $OH$ meet the sides of $ABC$ at $A'$, $B'$, $C'$, respectively. Let $AA_1$, $BB_1$, $CC_1$ be the altitudes of $ABC$. Consider the isogonal conjugation $A_2$ of $A_1$ wrt $AB'C'$. Similarly define $B_2$, $C_2$. Then the points $O$, $H$, $A_2$, $B_2$, $C_2$ lie on the same circle. 

\end{theorem}

\section{Magic triangle of facts}

Consider some geometric fact $\mathcal{F}$, so we can build a triangle $\mathcal{F}_1\mathcal{F}_2\mathcal{F}_3$ from three copies of this fact and look on the nice properties of the resulting configuration.

\bigskip

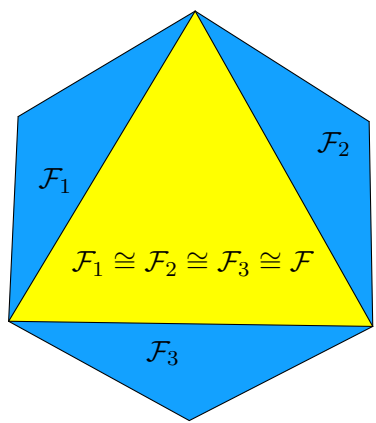
\definecolor{ffwwqq}{rgb}{0.07450980392156863,0.6313725490196078,1.0}
\definecolor{ffffqq}{rgb}{1.0,1.0,0.0}
\begin{tikzpicture}[line cap=round,line join=round,>=triangle 45,x=1.0cm,y=1.0cm]
\clip(-5.177577552691761,0.4874991744016925) rectangle (8.731983274655356,6.44317246840859);
\fill[line width=0.0pt,color=ffffqq,fill=ffffqq,fill opacity=1.0] (-0.48,1.92) -- (4.308979706863574,1.8518201093757962) -- (1.973535370739593,6.033288139039907) -- cycle;
\fill[line width=0.0pt,color=ffwwqq,fill=ffwwqq,fill opacity=1.0] (-0.48,1.92) -- (-0.35511084961342687,4.6339036313242215) -- (1.973535370739593,6.033288139039907) -- cycle;
\fill[line width=0.0pt,color=ffwwqq,fill=ffwwqq,fill opacity=1.0] (1.973535370739593,6.033288139039907) -- (4.261400283695918,4.56817914648597) -- (4.308979706863574,1.8518201093757962) -- cycle;
\fill[line width=0.0pt,color=ffwwqq,fill=ffwwqq,fill opacity=1.0] (-0.48,1.92) -- (4.308979706863574,1.8518201093757962) -- (1.8962256435206726,0.6030254706055125) -- cycle;
\draw (-0.48,1.92)-- (1.973535370739593,6.033288139039907);
\draw (4.308979706863574,1.8518201093757962)-- (1.973535370739593,6.033288139039907);
\draw (4.308979706863574,1.8518201093757962)-- (-0.48,1.92);
\draw (-0.35511084961342687,4.6339036313242215)-- (1.973535370739593,6.033288139039907);
\draw (-0.35511084961342687,4.6339036313242215)-- (-0.48,1.92);
\draw (-0.48,1.92)-- (1.8962256435206726,0.6030254706055125);
\draw (1.8962256435206726,0.6030254706055125)-- (4.308979706863574,1.8518201093757962);
\draw (4.308979706863574,1.8518201093757962)-- (4.261400283695918,4.56817914648597);
\draw (4.261400283695918,4.56817914648597)-- (1.973535370739593,6.033288139039907);
\draw (-0.22101137835150922,4.083115452135461) node[anchor=north west] {$\mathcal{F}_1$};
\draw (3.4440183410373555,4.569009543721105) node[anchor=north west] {$\mathcal{F}_2$};
\draw (1.1950228314123703,1.7924718775174235) node[anchor=north west] {$\mathcal{F}_3$};
\draw (0.214324663025794,2.9833386722187983) node[anchor=north west] {$\mathcal{F}_1\cong\mathcal{F}_2\cong\mathcal{F}_3\cong\mathcal{F}$};
\end{tikzpicture}

\begin{theorem}[\textbf{Feuerbach's Magic Triangle}]

Given a triangle $ABC$ and let $A'$ be the reflection of $A$ wrt midpoint of $BC$. Similarly define $B'$, $C'$. Let $\omega_A$, $\omega_B$, $\omega_C$ be the incircles of $A'BC$, $AB'C$, $ABC'$ and $\pi_A$, $\pi_B$, $\pi_C$ be its nine~- point circles. Consider the circle $\omega$ which is internally tangent to $\omega_A$, $\omega_B$, $\omega_C$. Consider the circle $\pi$ which is internally tangent to $\pi_A$, $\pi_B$, $\pi_C$. Then $\omega$ is tangent to $\pi$.

\end{theorem}

\begin{theorem}[\textbf{Magic Triangle related to the Simson lines}]

Consider any triangle $ABC$ and let $A'$, $B'$, $C'$ be a points outside of $ABC$ such that $A$, $B$, $C$ lie on the sides of $A'B'C'$. Consider any line $\ell$. Let $P_A$, $Q_A$ be the intersection points of $\ell$ with $(A'BC)$. Let the Simson lines of $P_A$, $Q_A$ wrt $A'BC$ meet at $X_A$. Similarly define $X_B$, $X_C$.

\begin{enumerate}
\item Points $X_A$, $X_B$, $X_C$ are collinear.
\item Let $(A'BC)$, $(AB'C)$, $(ABC')$ are concurrent at $P$ and $A_1B_1C_1$ be the pedal triangle of $P$ wrt $A'B'C'$. Then the midpoints of $A_1X_A$, $B_1X_B$, $C_1X_C$ are collinear.
\end{enumerate}
\end{theorem}

\begin{theorem}[\textbf{Magic Triangle related to the IMO 2011 G8}]

Consider any triangle $ABC$ and let $A'$ be the reflection of $A$ wrt midpoint of $BC$. Similarly define $B'$, $C'$. Let $H$ be the orthocenter of $ABC$, so $H$ lies on $(A'BC)$. Let $\ell_A$ be the tangent line through $H$ to $(A'BC)$. Let the reflections of $\ell_A$ wrt sides of $A'BC$ form a triangle with the circumcircle $\omega_A$. Consider the tangent point $T_A$ of $(A'BC)$ with $\omega_A$ (they are tangent because of IMO 2011 Problem G8). Similarly define $T_B$, $T_C$. Then $(T_AT_BT_C)$ is tangent to $(A'B'C')$.

\end{theorem}

\section{Incircles in some famous constructions}

In this section we will consider some famous configurations from plane geometry and add some incircles to them.

\begin{theorem}[\textbf{Mixtilinear circles}]

Given a triangle $ABC$ with the incenter $I$. Let $A$~-- mixtilinear incircle is tangent to $AB$, $AC$ at $A_B$, $A_C$, respectively. Similarly define the points $B_A$, $B_C$, $C_A$, $C_B$. Let $I_A$ be the incenter of $IC_BB_C$ and $W_A$ be the midpoint of smaller arc $BC$ of $(ABC)$. Similarly define the points $I_B$, $I_C$, $W_B$, $W_C$. Then $I_AW_A$, $I_BW_B$, $I_CW_C$ are concurrent.

\end{theorem}

\begin{theorem}[\textbf{Reim's theorem}]

Consider two arbitrary lines $l_1$, $l_2$ and let the points $A$, $B$, $C$ lie on $l_1$ and the points $D$, $E$, $F$ lie on $l_2$. Let also given that $ABDE$, $BCEF$ are cyclic. Consider the incenters $I_1$, $I_2$, $I_3$, $I_4$, $I_5$, $I_6$, $I_7$, $I_8$ of the triangles $ABD$, $AED$, $ABE$, $BED$, $BCE$, $BFE$, $BCF$, $CEF$, respectively. Let $P = I_1I_8\cap I_2I_7$, $Q = I_4I_5\cap I_3I_6$. Then the lines $l_1$, $l_2$, $PQ$ are concurrent.

\end{theorem}

\begin{theorem}[\textbf{Projectivity}]

Consider two arbitrary lines $l_1$, $l_2$ and let the points $A$, $B$, $C$ lie on $l_1$ and the points $D$, $E$, $F$ lie on $l_2$. Also let given that $AD$, $BE$, $CF$ are concurrent. Let $X = AE\cap DB$, $Y = AF\cap DC$, $Z = BF\cap CE$. Consider the incenters $I_1$, $I_2$, $I_3$, $I_4$ of $AXB$, $BZC$, $DXE$, $EZF$, respectively. Then in the case when $l_1$, $l_2$, $I_1I_2$ are concurrent at $P$ we get that $P$ lies on $I_3I_4$.

\end{theorem}

\begin{theorem}[\textbf{Feuerbach's theorem}]

Consider any triangle $ABC$ with the Feuerbach point $F$. Let $A'B'C'$, $A_1B_1C_1$ be the midpoint triangle and the orthic triangle of $ABC$, respectively. Let the incircle of $ABC$ is tangent to  its sides at $A_0$, $B_0$, $C_0$. Then

\begin{enumerate}
    \item The incenters of $FA'A_1$, $FB'B_1$, $FC'C_1$ lie on a circle which goes through $F$.
    \item The incenters of $FA'A_0$, $FB'B_0$, $FC'C_0$ lie on a circle which goes through $F$.
\end{enumerate}

\end{theorem}

\begin{theorem}[\textbf{Touching circles}]

Given a triangle $ABC$, let the circle $\omega_A$ through $B$, $C$ is tangent to the incircle of $ABC$ and meet $AB$, $AC$ second time at $A_B$, $A_C$, respectively. Consider the intersection $A' = BA_C\cap CA_B$ and let $I_A$ be the incenter of $A'BC$. Similarly define $I_B$, $I_C$. Then $AI_A$, $BI_B$, $CI_C$ are concurrent.

\end{theorem}

\begin{theorem}

Given a cyclic quadrilateral $ABCD$. Let $P = AB\cap CD$, $Q = BC\cap AD$. Let $l_P$ be the angle bisector of $\angle APD$ and $l_Q$ be the angle bisector of $\angle AQB$. Let $l_P$ meet $AD$, $BC$ at $X$, $Y$ and $l_Q$ meet $AB$, $CD$ at $Z$, $W$, respectively. Then the external homothety center of the incircles of the triangles $DXW$, $BYZ$ lies on $PQ$.

\end{theorem}

\bigskip

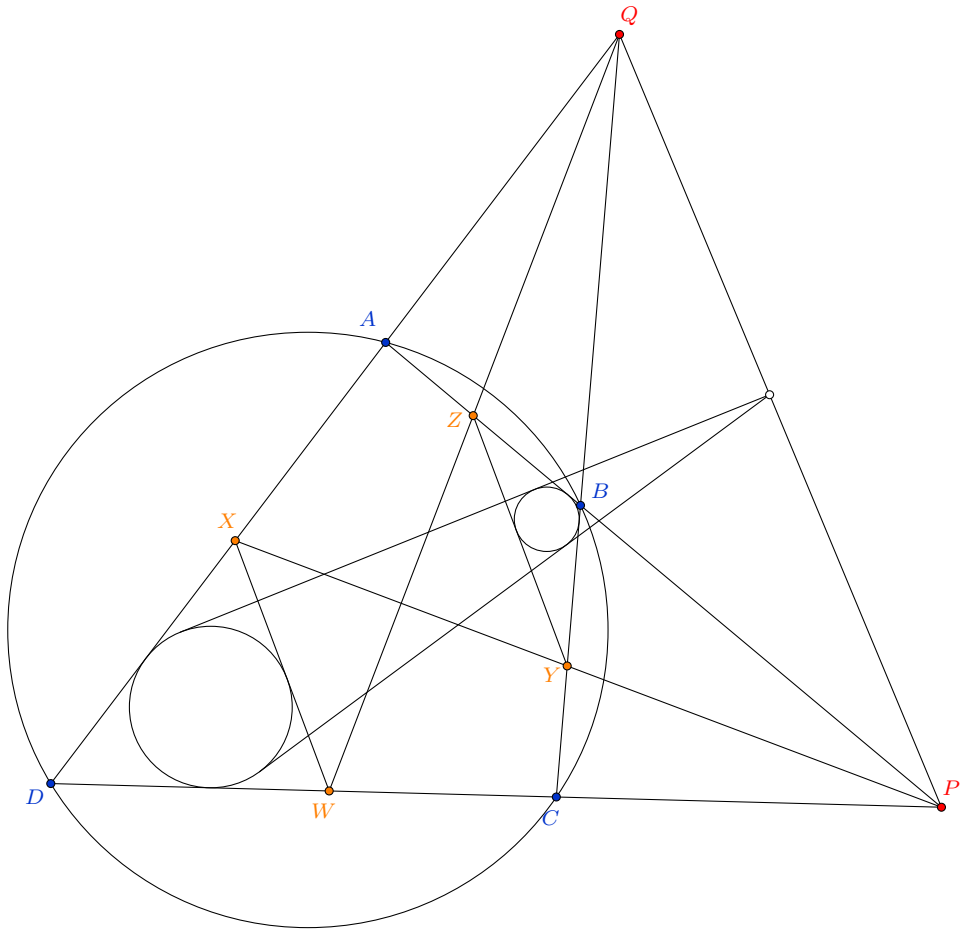
\definecolor{ffffff}{rgb}{1.0,1.0,1.0}
\definecolor{ffxfqq}{rgb}{1.0,0.4980392156862745,0.0}
\definecolor{ffqqqq}{rgb}{1.0,0.0,0.0}
\definecolor{qqqqff}{rgb}{0.0,0.2,0.8}
\begin{tikzpicture}[line cap=round,line join=round,>=triangle 45,x=1.0cm,y=1.0cm]
\draw(-8.054305160877124,8.42816678031475) circle (3.9474394042238194cm);
\draw(-4.912832089800837,9.895830717046344) circle (0.42699342010322244cm);
\draw(-9.331268606149544,7.405969007828541) circle (1.0707449239346156cm);
\draw (-9.010311707963254,9.6124306860057)-- (-7.774317585221203,6.293134695751864);
\draw (-5.880174972214633,11.270034787679263)-- (-4.6441808494725825,7.950738797425429);
\draw (-9.010311707963254,9.6124306860057)-- (0.27794508153676795,6.077442508510197);
\draw (-7.774317585221203,6.293134695751864)-- (-3.9562603180683684,16.32516155233317);
\draw (-11.435588123349302,6.391207434202445)-- (-3.9562603180683684,16.32516155233317);
\draw (-3.9562603180683684,16.32516155233317)-- (-4.786918891246252,6.21311264665423);
\draw (-3.9562603180683684,16.32516155233317)-- (0.27794508153676795,6.077442508510197);
\draw (0.27794508153676795,6.077442508510197)-- (-7.031417354104107,12.240774599744123);
\draw (-11.435588123349302,6.391207434202445)-- (0.27794508153676795,6.077442508510197);
\draw (-9.734277163749598,8.397976464772716)-- (-1.982130866741883,11.547329098759516);
\draw (-8.691451003039093,6.547407346163101)-- (-1.982130866741883,11.547329098759516);
\begin{scriptsize}
\draw [fill=qqqqff] (-4.786918891246252,6.21311264665423) circle (1.5pt);
\draw[color=qqqqff] (-4.859322067467473,5.939955760562543) node {$C$};
\draw [fill=qqqqff] (-11.435588123349302,6.391207434202445) circle (1.5pt);
\draw[color=qqqqff] (-11.647119838206953,6.21311264665423) node {$D$};
\draw [fill=qqqqff] (-7.031417354104107,12.240774599744123) circle (1.5pt);
\draw[color=qqqqff] (-7.266727676823075,12.550352403981368) node {$A$};
\draw [fill=qqqqff] (-4.4692451811702485,10.080323792200437) circle (1.5pt);
\draw[color=qqqqff] (-4.2076934814764835,10.274045019883976) node {$B$};
\draw [fill=ffqqqq] (-3.9562603180683684,16.32516155233317) circle (1.5pt);
\draw[color=ffqqqq] (-3.827576806315072,16.574863859065555) node {$Q$};
\draw [fill=ffqqqq] (0.27794508153676795,6.077442508510197) circle (1.5pt);
\draw[color=ffqqqq] (0.40800900262636364,6.340585860163684) node {$P$};
\draw [fill=ffxfqq] (-5.880174972214633,11.270034787679263) circle (1.5pt);
\draw[color=ffxfqq] (-6.126377651338843,11.220988891668492) node {$Z$};
\draw [fill=ffxfqq] (-4.6441808494725825,7.950738797425429) circle (1.5pt);
\draw[color=ffxfqq] (-4.8412212734121685,7.833843504131573) node {$Y$};
\draw [fill=ffxfqq] (-9.010311707963254,9.6124306860057) circle (1.5pt);
\draw[color=ffxfqq] (-9.113008670464213,9.873414920282835) node {$X$};
\draw [fill=ffxfqq] (-7.774317585221203,6.293134695751864) circle (1.5pt);
\draw[color=ffxfqq] (-7.845953086592845,6.031008055926439) node {$W$};
\draw [fill=ffffff] (-1.982130866741883,11.547329098759516) circle (1.5pt);
\end{scriptsize}
\end{tikzpicture}

\newpage

\begin{theorem}[\textbf{Isogonal conjugation}]

Consider two isogonal conjugated points $P$, $Q$ of a triangle $ABC$. Let $H_a$ be the external homothety center of the incircles of $BCP$, $BCQ$. Similarly define $H_b$, $H_c$. Then $H_a$, $H_b$, $H_c$ are collinear.

\end{theorem}

\begin{theorem}[\textbf{Arbelos}]

Consider any three collinear points $A$, $B$, $C$ ($B$ lies between $A$ and $C$) and let $\Omega_{AB}$, $\Omega_{BC}$, $\Omega_{CA}$ be the circles with diameters $AB$, $BC$, $CA$, respectively. Let the circle $\omega$ is tangent to $\Omega_{AB}$, $\Omega_{BC}$, $\Omega_{CA}$ at $C'$, $A'$, $B'$, respectively. Consider a triangle $A_1B_1C_1$ formed by the tangent lines from $A'$, $B'$, $C'$ to $\omega$. Then the common external tangent line to the $A_1$~-- excircle of $A_1B_1C_1$ and to $\Omega_{BC}$ is parallel to $AB$.

\end{theorem}

\begin{theorem}[\textbf{Confocal conics}]

Consider two ellipses $\mathcal{E}_1$ and $\mathcal{E}_2$ which share the same focuses $F_1$, $F_2$ and let $\mathcal{H}_1$, $\mathcal{H}_2$ be two hyperbolas with the same focuses $F_1$, $F_2$. Let $ABCD$ be a quadrilateral formed by the four intersection points of $\mathcal{E}_i$ and $\mathcal{H}_j$ ($A$, $B$ lie on $\mathcal{H}_1$ and $C$, $D$ lie on $\mathcal{H}_2$, respectively). Consider the incircles $\omega_1$, $\omega_2$ of the triangles $F_1AB$ and $F_2CD$. Then the external tangents to $\omega_1$, $\omega_2$ meet on $F_1F_2$.

\end{theorem}

\bigskip

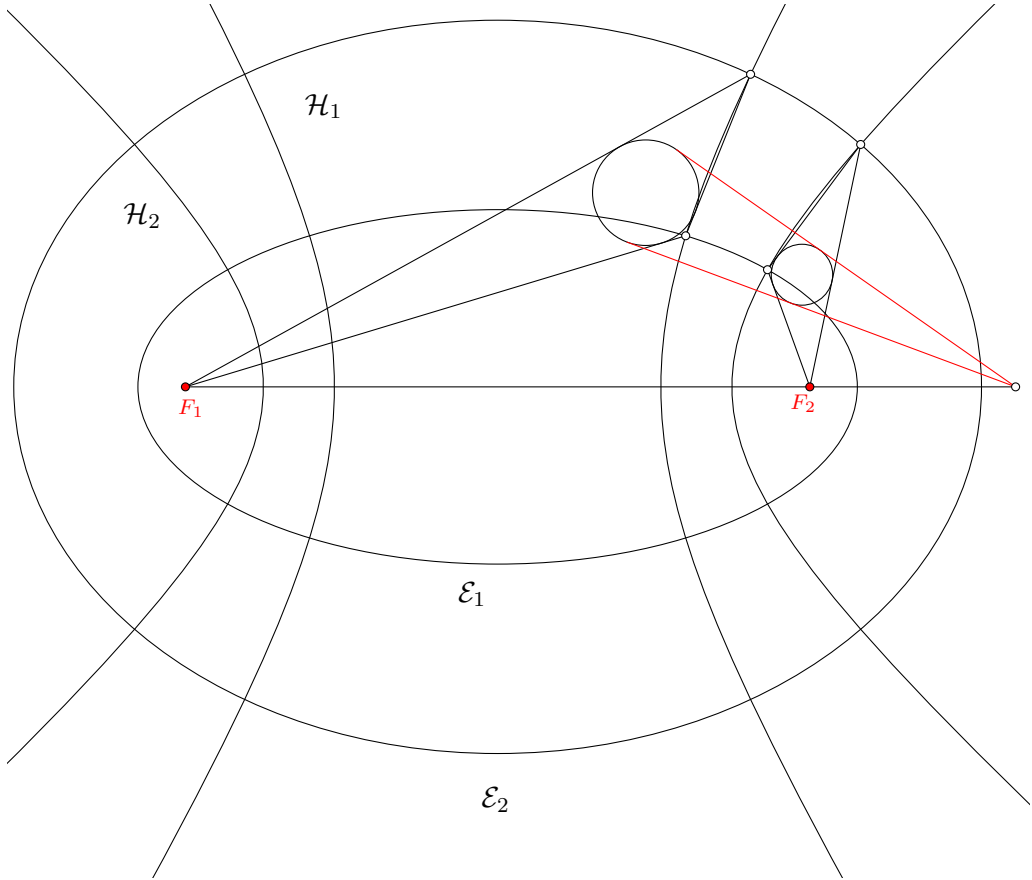
\definecolor{ffffff}{rgb}{1.0,1.0,1.0}
\definecolor{ffqqqq}{rgb}{1.0,0.0,0.0}
\begin{tikzpicture}[line cap=round,line join=round,>=triangle 45,x=1.0cm,y=1.0cm]
\clip(-3.0963860723379524,-6.509588141727433) rectangle (10.346078848221683,5.074471490275579);
\draw [rotate around={0.0:(3.35209222122451,0.0)}] (3.35209222122451,0.0) ellipse (6.363798673099702cm and 4.862056329411258cm);
\draw [samples=50,domain=-0.99:0.99,rotate around={0.0:(3.35209222122451,-0.0)},xshift=3.35209222122451cm,yshift=-0.0cm] plot ({3.0806139709565197*(1+(\x)^2)/(1-(\x)^2)},{2.7144353673139774*2*(\x)/(1-(\x)^2)});
\draw [samples=50,domain=-0.99:0.99,rotate around={0.0:(3.35209222122451,-0.0)},xshift=3.35209222122451cm,yshift=-0.0cm] plot ({3.0806139709565197*(-1-(\x)^2)/(1-(\x)^2)},{2.7144353673139774*(-2)*(\x)/(1-(\x)^2)});
\draw [rotate around={0.0:(3.3520922212245123,0.0)}] (3.3520922212245123,0.0) ellipse (4.730883965749434cm and 2.3500896357389487cm);
\draw [samples=50,domain=-0.99:0.99,rotate around={0.0:(3.35209222122451,-0.0)},xshift=3.35209222122451cm,yshift=-0.0cm] plot ({2.146385329943515*(1+(\x)^2)/(1-(\x)^2)},{3.5001959683395913*2*(\x)/(1-(\x)^2)});
\draw [samples=50,domain=-0.99:0.99,rotate around={0.0:(3.35209222122451,-0.0)},xshift=3.35209222122451cm,yshift=-0.0cm] plot ({2.146385329943515*(-1-(\x)^2)/(1-(\x)^2)},{3.5001959683395913*(-2)*(\x)/(1-(\x)^2)});
\draw(5.299794137326341,2.5747510091519428) circle (0.6996180107956123cm);
\draw(7.354672462836722,1.487102058870319) circle (0.404079221866037cm);
\draw (5.825197247990847,2.0034077985241643)-- (6.678815641564974,4.14481278461682);
\draw (6.678815641564974,4.14481278461682)-- (-0.7537988880692396,0.0);
\draw (-0.7537988880692396,0.0)-- (5.825197247990847,2.0034077985241643);
\draw (7.45798333051826,0.0)-- (6.9016328329746255,1.5536618614089621);
\draw (6.9016328329746255,1.5536618614089621)-- (8.126794369640141,3.214341858349191);
\draw (8.126794369640141,3.214341858349191)-- (7.45798333051826,0.0);
\draw (-0.7537988880692396,0.0)-- (10.164231363071245,6.618273626801106E-15);
\draw [color=ffqqqq] (10.164231363071245,6.618273626801106E-15)-- (5.053760472384541,1.9198213970981652);
\draw [color=ffqqqq] (5.703028797215133,3.146473992854693)-- (10.164231363071245,6.618273626801106E-15);
\draw (2.6956437436635494,-2.4520699283787843) node[anchor=north west] {$\mathcal{E}_1$};
\draw (3.0053779584229883,-5.177731018261846) node[anchor=north west] {$\mathcal{E}_2$};
\draw (0.6823713477271985,4.021375160093487) node[anchor=north west] {$\mathcal{H}_1$};
\draw (-1.7025821059204789,2.5656243507241245) node[anchor=north west] {$\mathcal{H}_2$};
\begin{scriptsize}
\draw [fill=ffqqqq] (-0.7537988880692396,0.0) circle (1.5pt);
\draw[color=ffqqqq] (-0.6804591972143313,-0.2839304250627125) node {$F_1$};
\draw [fill=ffqqqq] (7.45798333051826,0.0) circle (1.5pt);
\draw[color=ffqqqq] (7.372630386531074,-0.2219835821108247) node {$F_2$};
\draw [fill=ffffff] (8.126794369640141,3.214341858349191) circle (1.5pt);
\draw [fill=ffffff] (6.9016328329746255,1.5536618614089621) circle (1.5pt);
\draw [fill=ffffff] (5.825197247990847,2.0034077985241643) circle (1.5pt);
\draw [fill=ffffff] (6.678815641564974,4.14481278461682) circle (1.5pt);
\draw [fill=ffffff] (10.164231363071245,6.618273626801106E-15) circle (1.5pt);
\end{scriptsize}
\end{tikzpicture}

\newpage

\begin{theorem}[\textbf{Apollonian gasket}]

Let given three circles $\Omega_a$, $\Omega_b$, $\Omega_c$, which are pairwise tangent at $A'$, $B'$, $C'$. Let the circle $\Omega$ is internally tangent to $\Omega_a$, $\Omega_b$, $\Omega_c$ at $A$, $B$, $C$. Let $A_B$ be the second intersection point of $BC$ with $\Omega_b$, similarly define $B_A$, $C_A$, $A_C$, $C_B$, $B_C$. Let the line $l_A$ goes through the incenters of $A'A_BB$, $A'A_CC$, similarly define $l_B$, $l_C$. Then the triangle formed by $l_A$, $l_B$, $l_C$ is perspective to $ABC$.

\end{theorem}

\bigskip

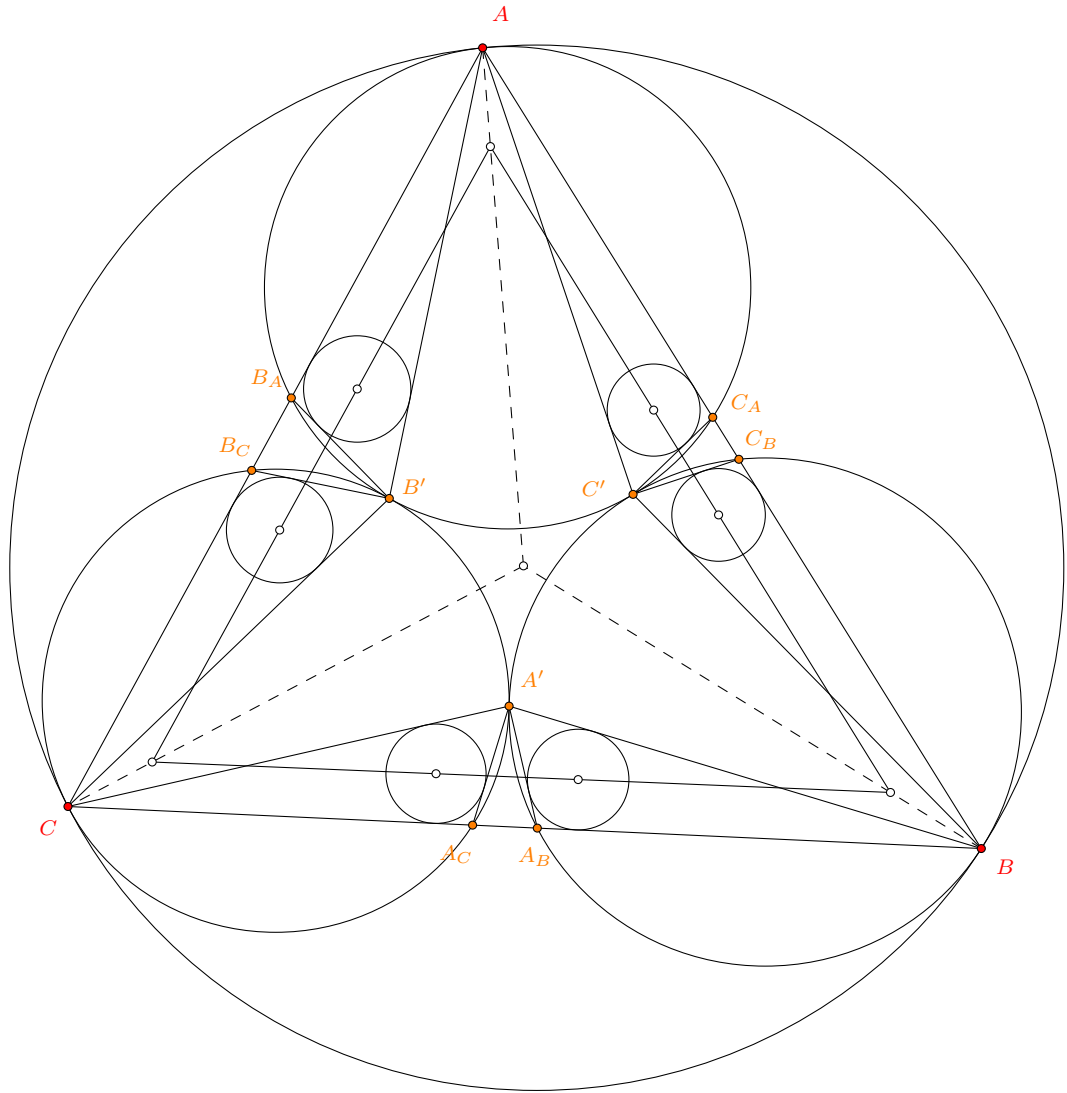
\definecolor{ffffff}{rgb}{1.0,1.0,1.0}
\definecolor{qqttcc}{rgb}{1.0,1.0,1.0}
\definecolor{ffqqqq}{rgb}{1.0,0.0,0.0}
\definecolor{ffxfqq}{rgb}{1.0,0.4980392156862745,0.0}
\begin{tikzpicture}[line cap=round,line join=round,>=triangle 45,x=1.0cm,y=1.0cm]
\draw(-13.16261847559686,-1.8757133431448623) circle (3.198204916250919cm);
\draw(-16.213225784891698,-7.352230729927838) circle (3.0706424772704226cm);
\draw(-9.774980392150134,-7.502138333912881) circle (3.3693478974809543cm);
\draw(-12.779122382671034,-5.588805419063611) circle (6.9310486067287265cm);
\draw (-14.71895996577586,-4.669691842038324)-- (-13.491188110828865,1.3055688953770597);
\draw (-13.491188110828865,1.3055688953770597)-- (-18.94477888525578,-8.754893757499119);
\draw (-18.94477888525578,-8.754893757499119)-- (-14.71895996577586,-4.669691842038324);
\draw (-18.94477888525578,-8.754893757499119)-- (-13.14341532996442,-7.423707957326654);
\draw (-13.14341532996442,-7.423707957326654)-- (-6.933079759338926,-9.312140004488322);
\draw (-6.933079759338926,-9.312140004488322)-- (-18.94477888525578,-8.754893757499119);
\draw (-11.512938437957121,-4.615616710007374)-- (-13.491188110828865,1.3055688953770597);
\draw (-13.491188110828865,1.3055688953770597)-- (-6.933079759338926,-9.312140004488322);
\draw (-6.933079759338926,-9.312140004488322)-- (-11.512938437957121,-4.615616710007374);
\draw(-15.140936040377113,-3.218157791883288) circle (0.7054921931886096cm);
\draw(-11.242138608404295,-3.4982369703518823) circle (0.6109190854748806cm);
\draw(-10.388298023135393,-4.8881761827036) circle (0.6148888057723509cm);
\draw(-12.234893570077144,-8.397807626369492) circle (0.6676527724006107cm);
\draw(-14.104447254220226,-8.320041708170116) circle (0.6586960669222557cm);
\draw(-16.160486737473978,-5.088758701919951) circle (0.7006283263125188cm);
\draw (-14.71895996577586,-4.669691842038324)-- (-16.007647267920674,-3.336646557627997);
\draw (-16.528690272097652,-4.297835964916392)-- (-14.71895996577586,-4.669691842038324);
\draw (-13.14341532996442,-7.423707957326654)-- (-13.623270407283057,-9.001768957800572);
\draw (-13.14341532996442,-7.423707957326654)-- (-12.772252920501264,-9.041249325703571);
\draw (-11.512938437957121,-4.615616710007374)-- (-10.121132511428009,-4.1506186838766865);
\draw (-10.46506961691447,-3.5937776147461524)-- (-11.512938437957121,-4.615616710007374);
\draw (-13.388885737550437,-0.0036174063685310604)-- (-8.127371891459351,-8.568664017159874);
\draw (-13.388885737550437,-0.0036174063685310604)-- (-17.83703920058495,-8.1647808871332);
\draw (-17.83703920058495,-8.1647808871332)-- (-8.127371891459351,-8.568664017159874);
\draw [dash pattern=on 4pt off 4pt] (-12.95440879684151,-5.563715779183114)-- (-13.491188110828865,1.3055688953770597);
\draw [dash pattern=on 4pt off 4pt] (-12.95440879684151,-5.563715779183114)-- (-18.94477888525578,-8.754893757499119);
\draw [dash pattern=on 4pt off 4pt] (-6.933079759338926,-9.312140004488322)-- (-12.95440879684151,-5.563715779183114);
\begin{scriptsize}
\draw [fill=ffxfqq] (-14.71895996577586,-4.669691842038324) circle (1.5pt);
\draw[color=ffxfqq] (-14.388000814685714,-4.504169557648878) node {$B'$};
\draw [fill=ffxfqq] (-11.512938437957121,-4.615616710007374) circle (1.5pt);
\draw[color=ffxfqq] (-12.022545700765956,-4.5337835364937956) node {$C'$};
\draw [fill=ffxfqq] (-13.14341532996442,-7.423707957326654) circle (1.5pt);
\draw[color=ffxfqq] (-12.832633068546695,-7.058466608103378) node {$A'$};
\draw [fill=ffqqqq] (-13.491188110828865,1.3055688953770597) circle (1.5pt);
\draw[color=ffqqqq] (-13.253878499792679,1.7580734692173348) node {$A$};
\draw [fill=ffqqqq] (-18.94477888525578,-8.754893757499119) circle (1.5pt);
\draw[color=ffqqqq] (-19.20217394732736,-9.03229026962206) node {$C$};
\draw [fill=ffqqqq] (-6.933079759338926,-9.312140004488322) circle (1.5pt);
\draw[color=ffqqqq] (-6.61674111572324,-9.55074618500173) node {$B$};
\draw [fill=ffxfqq] (-16.528690272097652,-4.297835964916392) circle (1.5pt);
\draw[color=ffxfqq] (-16.742297865140225,-3.9857136422692065) node {$B_C$};
\draw [fill=ffxfqq] (-16.007647267920674,-3.336646557627997) circle (1.5pt);
\draw[color=ffxfqq] (-16.323841949760553,-3.0784157903547833) node {$B_A$};
\draw [fill=ffxfqq] (-10.46506961691447,-3.5937776147461524) circle (1.5pt);
\draw[color=ffxfqq] (-10.013529028669721,-3.4052402533333886) node {$C_A$};
\draw [fill=ffxfqq] (-10.121132511428009,-4.1506186838766865) circle (1.5pt);
\draw[color=ffxfqq] (-9.819108060402344,-3.8885031581355182) node {$C_B$};
\draw [fill=ffxfqq] (-12.772252920501264,-9.041249325703571) circle (1.5pt);
\draw[color=ffxfqq] (-12.800229573835464,-9.421132206156813) node {$A_B$};
\draw [fill=ffxfqq] (-13.623270407283057,-9.001768957800572) circle (1.5pt);
\draw[color=ffxfqq] (-13.83714140459481,-9.388728711445584) node {$A_C$};
\draw [fill=qqttcc] (-15.140936040377113,-3.218157791883288) circle (1.5pt);
\draw [fill=qqttcc] (-11.242138608404295,-3.4982369703518823) circle (1.5pt);
\draw [fill=qqttcc] (-10.388298023135393,-4.8881761827036) circle (1.5pt);
\draw [fill=qqttcc] (-12.234893570077144,-8.397807626369492) circle (1.5pt);
\draw [fill=qqttcc] (-14.104447254220226,-8.320041708170116) circle (1.5pt);
\draw [fill=qqttcc] (-16.160486737473978,-5.088758701919951) circle (1.5pt);
\draw [fill=ffffff] (-17.83703920058495,-8.1647808871332) circle (1.5pt);
\draw [fill=ffffff] (-13.388885737550437,-0.0036174063685310604) circle (1.5pt);
\draw [fill=ffffff] (-8.127371891459351,-8.568664017159874) circle (1.5pt);
\draw [fill=ffffff] (-12.95440879684151,-5.563715779183114) circle (1.5pt);
\end{scriptsize}
\end{tikzpicture}

\newpage

\section{Looks like a crown}

Here we will create some crown looking-like facts.

\begin{theorem}[\textbf{Crown related to the orthocenters}]

Let $AA'$, $BB'$ be the altitudes of triangle $ABC$. Let the perpendicular from $A$ to $A'B'$ meet $BB'$ at $B_1$ and the perpendicular from $B$ to $A'B'$ meet $AA'$ at $A_1$. Consider the orthocenters $H_A$, $H_B$ of $ABA_1$, $ABB_1$, respectively. Let $H$, $H'$ be the orthocenters of $ABC$ and $CH_AH_B$. Then the lines $AB$, $H_AH_B$ meet on the perpendicular line from $C$ to $HH'$.
\end{theorem}

\begin{theorem}[\textbf{Crown related to the incircles}]

Given a triangle $ABC$, let its incircle is tangent to its sides at $A'$, $B'$, $C'$, respectively. Let the circle $\omega_A$ with center at $I_A$ is tangent to $A'B'$ at $T_A$ and also is tangent to $AA'$, $AB$ (see the picture below). Let the second tangent line from $B$ to $\omega_A$ meet $AA'$ at $A_1$. Similarly define $\omega_B$, $T_B$, $I_B$ and $B_1$. Let $I$ be the incenter of $CA'B'$. Then $A_1T_A$, $B_1T_B$, $I_AI_B$, $IC'$ are concurrent.

\end{theorem}

\begin{theorem}[\textbf{Crown related to the tangent circles}]

Given a triangle $ABC$ with the incircle $\omega$. Let the circle $\Omega$ goes through $B$, $C$ and is tangent to $\omega$. Let $\Omega$ meet $CA$, $CB$ second time at $B'$, $A'$, respectively. Consider the circle $\omega_A$ which is tangent to $AA'$, $AB$ and is internally tangent to $\Omega$. Let the second tangent line from $B$ to $\omega_A$ meet $AA'$ at $A_1$ and also meet $\Omega$ second time at $A_2$. Similarly define $B_1$, $B_2$. Then $A'B'$, $A_1B_1$, $A_2B_2$ are concurrent.

\end{theorem}

\begin{theorem}[\textbf{Crown related to the isogonal conjugation}]

Given a triangle $ABC$ and let $P$ be any point on its Neuberg cubic. Let $A_1$ be the isogonal conjugation of $A$ wrt $BPC$. Consider the intersection point $A'$ of $BC$, $PA_1$. Similarly define $B_1$, $B'$. Consider the intersection points $A_2$, $B_2$ of $A'B'$ with $A_1B$, $B_1A$, respectively. Let $A_1'$, $B_1'$, $A_2'$, $B_2'$ be the isogonal conjugations of $A_1$, $B_1$, $A_2$, $B_2$ wrt $ABC$. Let $P_C = A_1'B_1'\cap A_2'B_2'$, similarly define $P_A$, $P_B$. Then $P_A$, $P_B$, $P_C$ are collinear.

\end{theorem}

\bigskip

\definecolor{ffffqq}{rgb}{1.0,1.0,0.0}
\definecolor{ffffff}{rgb}{1.0,1.0,1.0}
\definecolor{ffqqqq}{rgb}{1.0,0.0,0.0}
\begin{tikzpicture}[line cap=round,line join=round,>=triangle 45,x=1.0cm,y=1.0cm]
\fill[line width=0.0pt,color=ffffqq,fill=ffffqq,fill opacity=0.7] (-24.576887031417943,28.73328833156983) -- (-23.718445487048896,32.04515327206028) -- (-22.650191121461113,31.21620993659309) -- (-21.240037739919654,33.033466245350354) -- (-20.20369988770961,30.58207486444369) -- (-18.762676113683664,31.191274804767108) -- (-19.4082763508622,28.700550998225975) -- cycle;
\draw (-24.576887031417943,28.73328833156983)-- (-21.240037739919654,33.033466245350354);
\draw (-21.240037739919654,33.033466245350354)-- (-19.4082763508622,28.700550998225975);
\draw (-24.576887031417943,28.73328833156983)-- (-19.4082763508622,28.700550998225975);
\draw (-18.78774902213501,27.23273300623147)-- (-19.4082763508622,28.700550998225975);
\draw (-23.73253042601416,29.821405686813566)-- (-18.78774902213501,27.23273300623147);
\draw (-24.576887031417943,28.73328833156983)-- (-23.718445487048896,32.04515327206028);
\draw (-19.4082763508622,28.700550998225975)-- (-18.762676113683664,31.191274804767108);
\draw (-24.576887031417943,28.73328833156983)-- (-18.762676113683664,31.191274804767108);
\draw (-19.4082763508622,28.700550998225975)-- (-23.718445487048896,32.04515327206028);
\draw (-23.852529322406493,31.527858231623124)-- (-19.00136168676419,30.27042656941373);
\draw (-22.691118996823374,30.261667528394906)-- (-21.258388662063116,30.136199959409435);
\draw [dash pattern=on 3pt off 3pt] (-21.61825551627652,28.71454872974274)-- (-21.240037739919654,33.033466245350354);
\draw (-24.62851372607856,28.760599267270628) node[anchor=north west] {$A$};
\draw (-19.350089239280514,28.800226535751104) node[anchor=north west] {$B$};
\draw (-21.19635430641167,33.70677879078023) node[anchor=north west] {$C$};
\draw (-21.259474479646922,30.234733654744396) node[anchor=north west] {$H$};
\draw (-22.687568399094523,30.33776455279364) node[anchor=north west] {$H'$};
\draw (-23.752721322439424,29.870162784724002) node[anchor=north west] {$H_B$};
\draw (-18.78774902213501,27.9) node[anchor=north west] {$H_A$};
\draw (-23.73694127913061,32.58253160638629) node[anchor=north west] {$B_1$};
\draw (-18.876687940016115,31.702806246119682) node[anchor=north west] {$A_1$};
\draw (-17.038312894539366,31.169937190883672) node[anchor=north west] {$\textbf{ORTHOCENTERS}$};
\begin{scriptsize}
\draw [fill=ffqqqq] (-24.576887031417943,28.73328833156983) circle (1.5pt);
\draw [fill=ffqqqq] (-21.240037739919654,33.033466245350354) circle (1.5pt);
\draw [fill=ffqqqq] (-19.4082763508622,28.700550998225975) circle (1.5pt);
\draw [fill=ffffff] (-22.650191121461113,31.21620993659309) circle (1.5pt);
\draw [fill=ffffff] (-20.20369988770961,30.58207486444369) circle (1.5pt);
\draw [fill=ffffff] (-23.718445487048896,32.04515327206028) circle (1.5pt);
\draw [fill=ffffff] (-18.762676113683664,31.191274804767108) circle (1.5pt);
\draw [fill=ffffff] (-19.00136168676419,30.27042656941373) circle (1.5pt);
\draw [fill=ffffff] (-23.852529322406493,31.527858231623124) circle (1.5pt);
\draw [fill=ffffff] (-18.78774902213501,27.23273300623147) circle (1.5pt);
\draw [fill=ffffff] (-23.73253042601416,29.821405686813566) circle (1.5pt);
\draw [fill=ffffff] (-22.691118996823374,30.261667528394906) circle (1.5pt);
\draw [fill=ffffff] (-21.61825551627652,28.71454872974274) circle (1.5pt);
\draw [fill=ffffff] (-21.258388662063105,30.13619995940944) circle (1.5pt);
\end{scriptsize}
\end{tikzpicture}

\newpage

\definecolor{ffffqq}{rgb}{1.0,1.0,0.0}
\definecolor{ffffff}{rgb}{1.0,1.0,1.0}
\definecolor{ffqqqq}{rgb}{1.0,0.0,0.0}
\begin{tikzpicture}[line cap=round,line join=round,>=triangle 45,x=1.0cm,y=1.0cm]
\fill[line width=0.0pt,color=ffffqq,fill=ffffqq,fill opacity=0.7] (-20.207705220699907,22.055748319428762) -- (-21.04162259210748,25.472505167841927) -- (-19.074289963428726,24.624502148801493) -- (-17.970883479375622,27.125244584623413) -- (-16.14284536467463,25.09312863060617) -- (-14.713852343011744,26.160913356904928) -- (-13.597718374266407,22.263869699615885) -- cycle;
\draw (-20.207705220699907,22.055748319428762)-- (-17.970883479375622,27.125244584623413);
\draw (-17.970883479375622,27.125244584623413)-- (-13.597718374266407,22.263869699615885);
\draw(-17.457034960978614,23.91091820248667) circle (1.7676865648038742cm);
\draw(-19.235318136569415,23.334948651624426) circle (1.247965410973432cm);
\draw(-15.547847666705536,23.686308881015226) circle (1.4831057907684964cm);
\draw(-17.736077708533152,25.656441315030687) circle (0.8077536308783734cm);
\draw [dash pattern=on 3pt off 3pt] (-19.235318136569415,23.334948651624426)-- (-15.547847666705536,23.686308881015226);
\draw [dash pattern=on 3pt off 3pt] (-17.401405327531915,22.14410719525621)-- (-17.736077708533152,25.656441315030687);
\draw [dash pattern=on 3pt off 3pt] (-17.530354254922017,23.49740595276201)-- (-21.04162259210748,25.472505167841927);
\draw [dash pattern=on 3pt off 3pt] (-17.530354254922017,23.49740595276201)-- (-14.713852343011744,26.160913356904928);
\draw (-20.207705220699907,22.055748319428762)-- (-13.597718374266407,22.263869699615885);
\draw (-13.597718374266407,22.263869699615885)-- (-14.713852343011744,26.160913356904928);
\draw (-20.207705220699907,22.055748319428762)-- (-21.04162259210748,25.472505167841927);
\draw (-20.76854818936889,24.353654706544376)-- (-14.483984595216334,25.358317366059477);
\draw (-13.597718374266407,22.263869699615885)-- (-21.04162259210748,25.472505167841927);
\draw (-20.207705220699907,22.055748319428762)-- (-14.713852343011744,26.160913356904928);
\draw (-20.2096434171974,22.077577651617582) node[anchor=north west] {$A$};
\draw (-13.596247978859644,22.386152046677783) node[anchor=north west] {$B$};
\draw (-17.91520010185573,27.5) node[anchor=north west] {$C$};
\draw (-21.092121615405738,25.890201883689187) node[anchor=north west] {$B_1$};
\draw (-14.831717456351313,26.636782508375815) node[anchor=north west] {$A_1$};
\draw (-17.728322365764555,26.01534652072531) node[anchor=north west] {$I$};
\draw (-19.66977440182289,24.61789807382195) node[anchor=north west] {$T_B$};
\draw (-15.786870329706217,25.27490741826159) node[anchor=north west] {$T_A$};
\draw (-17.43762366517828,22.18800637137059) node[anchor=north west] {$C'$};
\draw (-19.347929411888085,23.429024021978798) node[anchor=north west] {$I_B$};
\draw (-15.54808211136749,23.814886652840173) node[anchor=north west] {$I_A$};
\draw (-12.298485922670917,25.20190637999052) node[anchor=north west] {$\textbf{INCIRCLES}$};
\begin{scriptsize}
\draw [fill=ffqqqq] (-20.207705220699907,22.055748319428762) circle (1.5pt);
\draw [fill=ffqqqq] (-17.970883479375622,27.125244584623413) circle (1.5pt);
\draw [fill=ffqqqq] (-13.597718374266407,22.263869699615885) circle (1.5pt);
\draw [fill=ffffff] (-17.401405327531915,22.14410719525621) circle (1.5pt);
\draw [fill=ffffff] (-19.074289963428726,24.624502148801493) circle (1.5pt);
\draw [fill=ffffff] (-16.14284536467463,25.09312863060617) circle (1.5pt);
\draw [fill=ffffff] (-15.547847666705536,23.686308881015226) circle (1.5pt);
\draw [fill=ffffff] (-20.76854818936889,24.353654706544376) circle (1.5pt);
\draw [fill=ffffff] (-21.04162259210748,25.472505167841927) circle (1.5pt);
\draw [fill=ffffff] (-14.713852343011744,26.160913356904928) circle (1.5pt);
\draw [fill=ffffff] (-14.483984595216334,25.358317366059477) circle (1.5pt);
\draw [fill=ffffff] (-17.736077708533152,25.656441315030687) circle (1.5pt);
\draw [fill=ffffff] (-19.235318136569415,23.334948651624426) circle (1.5pt);
\draw [fill=ffffff] (-19.432318953379767,24.56726693241154) circle (1.5pt);
\draw [fill=ffffff] (-15.78196717805454,25.1508193251626) circle (1.5pt);
\draw [fill=ffffff] (-17.530354254922017,23.49740595276201) circle (1.5pt);
\end{scriptsize}
\end{tikzpicture}

\definecolor{ffffqq}{rgb}{1.0,1.0,0.0}
\definecolor{ffffff}{rgb}{1.0,1.0,1.0}
\definecolor{ffqqqq}{rgb}{1.0,0.0,0.0}
\begin{tikzpicture}[line cap=round,line join=round,>=triangle 45,x=1.0cm,y=1.0cm]
\fill[line width=0.0pt,color=ffffqq,fill=ffffqq,fill opacity=0.7] (-25.28331908758691,17.18514207937873) -- (-24.20664813785236,21.733952164601696) -- (-23.020545062812786,20.473788782421877) -- (-21.526143276911135,22.64570672580462) -- (-20.624491502880225,19.95376250782244) -- (-19.6005991504498,20.562235182958947) -- (-19.571189191327502,16.809055039776826) -- cycle;
\draw (-25.28331908758691,17.18514207937873)-- (-21.526143276911135,22.64570672580462);
\draw (-21.526143276911135,22.64570672580462)-- (-19.571189191327502,16.809055039776826);
\draw (-25.28331908758691,17.18514207937873)-- (-19.571189191327502,16.809055039776826);
\draw(-22.075563813752634,18.739387366354578) circle (1.76162999377408cm);
\draw(-22.386584189640686,17.614806727538607) circle (2.928427014432468cm);
\draw(-20.914261078263593,18.232245518910705) circle (1.3318788612638288cm);
\draw(-23.368380387024718,18.599967660442257) circle (1.5375762543157205cm);
\draw (-25.28331908758691,17.18514207937873)-- (-24.20664813785236,21.733952164601696);
\draw (-24.20664813785236,21.733952164601696)-- (-19.571189191327502,16.809055039776826);
\draw (-19.6005991504498,20.562235182958947)-- (-25.28331908758691,17.18514207937873);
\draw (-19.6005991504498,20.562235182958947)-- (-19.571189191327502,16.809055039776826);
\draw (-25.772465949980494,17.458064494553508) node[anchor=north west] {$A$};
\draw (-19.472724145783854,17.05030425251753) node[anchor=north west] {$B$};
\draw (-21.587871041109484,23.127210958150153) node[anchor=north west] {$C$};
\draw (-24.17601922131145,22.072193549171153) node[anchor=north west] {$B_1$};
\draw (-19.579549746557877,20.816721225007747) node[anchor=north west] {$A_1$};
\draw (-23.136842252332798,20.92402655185932) node[anchor=north west] {$B'$};
\draw (-20.658488314375496,20.40896098297177) node[anchor=north west] {$A'$};
\draw (-25.591750746497397,19.593440498899813) node[anchor=north west] {$B_2$};
\draw (-19.558184626403072,18.72426735140207) node[anchor=north west] {$A_2$};
\draw (-17.582495641943987,20.55918844056397) node[anchor=north west] {$\textbf{TANGENT CIRCLES}$};
\begin{scriptsize}
\draw [fill=ffqqqq] (-25.28331908758691,17.18514207937873) circle (1.5pt);
\draw [fill=ffqqqq] (-21.526143276911135,22.64570672580462) circle (1.5pt);
\draw [fill=ffqqqq] (-19.571189191327502,16.809055039776826) circle (1.5pt);
\draw [fill=ffffff] (-23.020545062812786,20.473788782421877) circle (1.5pt);
\draw [fill=ffffff] (-20.624491502880225,19.95376250782244) circle (1.5pt);
\draw [fill=ffffff] (-24.20664813785236,21.733952164601696) circle (1.5pt);
\draw [fill=ffffff] (-19.6005991504498,20.562235182958947) circle (1.5pt);
\draw [fill=ffffff] (-24.783360822725147,19.297407853758962) circle (1.5pt);
\draw [fill=ffffff] (-19.584161899924705,18.464579690724076) circle (1.5pt);
\end{scriptsize}
\end{tikzpicture}

\definecolor{ffffqq}{rgb}{1.0,1.0,0.0}
\definecolor{ffffff}{rgb}{1.0,1.0,1.0}
\definecolor{ffqqqq}{rgb}{1.0,0.0,0.0}
\begin{tikzpicture}[line cap=round,line join=round,>=triangle 45,x=1.0cm,y=1.0cm]
\fill[line width=0.pt,color=ffffqq,fill=ffffqq,fill opacity=0.7] (-0.2972929668780835,1.32889990133897) -- (-0.13481981603201956,4.547285591360543) -- (0.9362777809740461,3.679483531104444) -- (2.008948273431996,5.72346982381959) -- (3.4751176645838777,3.9404899612127804) -- (4.347255732262364,5.00829907632797) -- (5.622664862579285,1.32889990133897) -- cycle;
\draw (-0.2972929668780835,1.32889990133897)-- (2.008948273431996,5.72346982381959);
\draw (2.008948273431996,5.72346982381959)-- (5.622664862579285,1.32889990133897);
\draw (-0.2972929668780835,1.32889990133897)-- (5.622664862579285,1.32889990133897);
\draw (2.1360639323467243,6.32272935870331) node[anchor=north west] {$C$};
\draw (-0.769436842847076,2.000796955602536) node[anchor=north west] {$A$};
\draw (5.731621141649053,2.000796955602536) node[anchor=north west] {$B$};
\draw (2.3358171106412984,2.545578350951373)-- (-0.13481981603201956,4.547285591360543);
\draw (2.3358171106412984,2.545578350951373)-- (4.347255732262364,5.00829907632797);
\draw (5.622664862579285,1.32889990133897)-- (4.347255732262364,5.00829907632797);
\draw (-0.2972929668780835,1.32889990133897)-- (-0.13481981603201956,4.547285591360543);
\draw (-0.18444534847673497,3.564267150171479)-- (4.674649447007495,4.063808296627741);
\draw (-0.07938040873854844,5.124210288935869) node[anchor=north west] {$B_1$};
\draw (-1.11569916842847094,3.744097420718815) node[anchor=north west] {$B_2$};
\draw (4.297030133897113,5.632672924594783) node[anchor=north west] {$A_1$};
\draw (4.714695870331222,4.270719436222691) node[anchor=north west] {$A_2$};
\draw (0.8285885835095143,3.7077786610288923) node[anchor=north west] {$B'$};
\draw (3.3527423819591284,3.9620099788583496) node[anchor=north west] {$A'$};
\draw (2.3358171106412984,2.763490909090908) node[anchor=north west] {$P$};
\draw (7.529399746300217,4.161763157152923) node[anchor=north west] {$\textbf{ISOGONAL POINTS}$};
\begin{scriptsize}
\draw [fill=ffqqqq] (-0.2972929668780835,1.32889990133897) circle (1.5pt);
\draw [fill=ffqqqq] (2.008948273431996,5.72346982381959) circle (1.5pt);
\draw [fill=ffqqqq] (5.622664862579285,1.32889990133897) circle (1.5pt);
\draw [fill=ffffff] (2.3358171106412984,2.545578350951373) circle (1.5pt);
\draw [fill=ffffff] (4.347255732262364,5.00829907632797) circle (1.5pt);
\draw [fill=ffffff] (-0.13481981603201956,4.547285591360543) circle (1.5pt);
\draw [fill=ffffff] (0.9362777809740461,3.679483531104444) circle (1.5pt);
\draw [fill=ffffff] (3.4751176645838777,3.9404899612127804) circle (1.5pt);
\draw [fill=ffffff] (-0.18444534847673497,3.564267150171479) circle (1.5pt);
\draw [fill=ffffff] (4.674649447007495,4.063808296627741) circle (1.5pt);
\end{scriptsize}
\end{tikzpicture}

\end{document}